\documentclass[lettersize,journal]{IEEEtran}
\usepackage{amsmath,amsfonts}
\usepackage{algorithmic}
\usepackage{algorithm}
\usepackage{array}
\usepackage[caption=false,font=normalsize,labelfont=sf,textfont=sf]{subfig}
\usepackage{textcomp}
\usepackage{stfloats}
\usepackage{url}
\usepackage{verbatim}
\usepackage{graphicx}
\usepackage{cite}
\usepackage{color,multirow}
\usepackage{stfloats}
\newtheorem{theorem}{Theorem}[section]
\newtheorem{lemma}[theorem]{Lemma}
\newtheorem{definition}[theorem]{Definition}
\UseRawInputEncoding

\newcommand{\rmF}{F}
\newcommand{\rmH}{H}
\newcommand{\rmT}{T}

\begin{document}

\title{Efficient Orthogonal Decomposition with Automatic Basis Extraction for Low-Rank Matrix Approximation}

\author{Weijie Shen, Weiwei Xu and Lei Zhu,~\IEEEmembership{Member,~IEEE}
\thanks{The work was supported in part by the major key project of Peng Cheng Laboratory under grant PCL2023AS1-2.}
\thanks{Weijie Shen is with the School of Mathematics and Statistics, Nanjing University of Information Science and Technology, Nanjing 210044, China (e-mail:swj@nuist.edu.cn).}
\thanks{Weiwei Xu is with the School of Mathematics and Statistics, Nanjing University of Information Science and Technology, Nanjing 210044, China. The Peng Cheng Laboratory, Shenzhen 518055, China, and with the Pazhou Laboratory (Huangpu), Guangzhou 510555, China (e-mail:wwxu@nuist.edu.cn).}
\thanks{Lei Zhu is with the College of Engineering, Nanjing Agricultural University, Nanjing 210031, China (e-mail:zhulei@njau.edu.cn).}
}

\markboth{Journal of \LaTeX\ Class Files,~Vol.~14, No.~8, August~2023}%
{Shell \MakeLowercase{\textit{et al.}}: A Sample Article Using IEEEtran.cls for IEEE Journals}


\maketitle

\begin{abstract}
Low-rank matrix approximation play a ubiquitous role in various applications such as image processing, signal processing, and data analysis. Recently, random algorithms of low-rank matrix approximation have gained widespread adoption due to their speed, accuracy, and robustness, particularly in their improved implementation on modern computer architectures.
Existing low-rank approximation algorithms often require prior knowledge of the rank of the matrix, which is typically unknown. To address this bottleneck,
we propose a low-rank approximation algorithm termed efficient orthogonal decomposition with automatic basis extraction (EOD-ABE) tailored for the scenario where the rank of the matrix is unknown.
Notably, we introduce a randomized algorithm to automatically extract the basis that reveals the rank.
The efficacy of the proposed algorithms is theoretically and numerically validated, demonstrating superior speed, accuracy, and robustness compared to existing methods. Furthermore, we apply the algorithms to image reconstruction, achieving remarkable results.
\end{abstract}


\section{Introduction}
Low-rank matrix approximation, that is, given a matrix $\mathbf{A}$, and the accuracy parameter $\varepsilon$, we seek a low-rank approximation $\tilde{\mathbf{A}}$ of $\mathbf{A}$, such as:
\begin{eqnarray*}
  \|\mathbf{A}-\tilde{\mathbf{A}}\|\leq(1+\varepsilon)\|\mathbf{A}-\mathbf{A}_{d}\|,
\end{eqnarray*}
where $d$ is a parameter, and $A_{d}$ is a best rank-$d$ approximation of $\mathbf{A}$.
Such compact representation that preserves the most important information of high-dimensional matrices can significantly reduce memory requirements and, more importantly, computational costs when the latter scales with the dimensionality.

Singular value decomposition (SVD), CPQR, UTV, and p-QLP all provide complete decompositions of the input matrix, see \cite{kc20,kc21,kl18,hmt11}. While these techniques are effective for achieving precise matrix decompositions, they come with high arithmetic costs due to the extensive computations required for complete decompositions, including eigenvalue and eigenvector calculations, pivoting, and iterative processes. The communication costs are also substantial, as handling large matrices in parallel or distributed computing environments necessitates significant data transfer and synchronization across nodes. These costs can limit the practicality of these methods for large-scale applications.
Recently developed algorithms for low-rank approximation \cite{s19} \cite{di19} \cite{kc19} based on random sampling schemes have been shown to be remarkably computationally efficient, highly accurate and robust, and are known to outperform the traditional algorithms in many practical situations.
The strategy of the random method is as follows: first, a random sampling method is applied to transform the input matrix into a low-dimensional space. Second, complete decomposition is performed on the reduced-size matrix, which is a key step in the process of this method. Finally, the processed matrix is projected back to the original space.
Compared with classical methods, random methods can better exploit parallel architectures in the case of low-rank approximation.

A bottleneck associated with existing randomized methods, however, is their requirement for prior knowledge of the matrix rank. This can be challenging in practical applications, especially for large matrices, as it requires significant time to calculate the rank in advance.
Therefore, this study aims to design a new random matrix decomposition method that can automatically extract the basis of the matrix and then determine the rank of the matrix, thereby overcoming this bottleneck.

\subsection{Applications}
Low-rank matrix approximation plays a significant role in various applications, including image processing, signal processing, and data analysis. Matrices with low-rank structure are widely used in image reconstruction \cite{kc20} \cite{kc21}, genomics \cite{dgme17} \cite{us19}, sparse matrix problems \cite{drs16}, mechanical fault diagnosis \cite{mhl23}, image inpainting \cite{hwz20}, principal component analysis \cite{kl18} \cite{kclw21}, video noise-reduction \cite{lsw20}, image classification \cite{z24}, image denoising \cite{xls20}, neural calcium imaging video segmentation and hyperspectral compressed recovery \cite{hv20}, etc.

\subsection{Contributions}
We propose a low-rank approximation algorithm called efficient orthogonal decomposition with automatic basis extraction (EOD-ABE), designed for matrix with unknown numerical rank.
\begin{itemize}
  \item EOD-ABE uses random sampling to automatically extract bases from matrices with low numerical rank, to obtain approximations.
Given a matrix $\mathbf{A}\in\mathbb{C}^{m\times n}$ with unknown rank. The randomized algorithm for basis extraction builds an orthogonal matrix $\mathbf{Q}$ to make
\begin{eqnarray*}
  \|(\mathbf{I}-\mathbf{Q}\mathbf{Q}^{H})\mathbf{A}\|\leq\varepsilon. 
\end{eqnarray*}
As a result, the matrix decomposition of the $m\times n$ matrix $\mathbf{A}$ is converted into the matrix decomposition of the $r\times n$ matrix $\mathbf{Q}^{H}\mathbf{A}$, so as to efficiently calculate the computer.
  \item EOD-ABE generates an approximation $\tilde{\mathbf{A}}$ such as:
\begin{equation}\label{EQLP}
  \tilde{\mathbf{A}} = \mathbf{U}\mathbf{D}\mathbf{V}^{H},
\end{equation}
where $\mathbf{U}\in\mathbb{C}^{m\times r}$ and $\mathbf{V}\in\mathbb{C}^{n\times r}$ are column orthonormal matrices that constitute approximations to the numerical range of $\mathbf{A}$ and $\mathbf{A}^{H}$, respectively. $\mathbf{D}$ is an upper triangular matrix, and its diagonals constitute approximations to the first $r$ singular values of $\mathbf{A}$.
  \item Through theoretical analysis and numerical experiments, we show that EOD-ABE can perform low-rank approximation of matrices of unknown rank while determining the rank. Because EOD-ABE uses a random method to automatically extract bases, it has low computational cost, high accuracy, and is compatible with modern computing environments.
\end{itemize}

\subsection{Notation}
Throughout this paper, by
$\mathbb{R}^{m\times n}$ and $\mathbb{C}^{m\times n}$ we denote the sets of $m\times n$ matrices with entries in real number field, and $m\times n$ matrices with entries in complex number field, respectively.
The symbols $\mathbf{I}_{n}$ and $\mathbf{O}_{m\times n}$ represent the identity matrix of order $n$ and $m\times n$ zero matrix, respectively.
For a real matrix $\mathbf{A}=(a_{ij})\in\mathbb{R}^{m\times n}$, by $\mathbf{A}^{\rmT}$ we denote the transpose.
For a complex matrix $\mathbf{A}=(a_{ij})\in\mathbb{C}^{m\times n}$, by $\mathbf{A}^{\rmH}$, $\mathbf{A}^{-1}$, $rank(\mathbf{A})$ and $\mathrm{tr}(\mathbf{A})$ we denote the conjugate transpose, inverse, rank and trace of matrix $\mathbf{A}$, respectively. We use $\|\cdot\|,\|\cdot\|_{2}$ and $\|\cdot\|_{\rmF}$ to denote the unitary invariant norm, spectral norm and Frobenius norm of a matrix, respectively. The singular value set of $\mathbf{A}$ is denoted by $\sigma(\mathbf{A})$. $\mathbf{A}\preceq \mathbf{B}$ denotes $\mathbf{A}-\mathbf{B}$ is negative semi-definite. $\mathbf{A}\succeq \mathbf{B}$ denotes $\mathbf{A}-\mathbf{B}$ is positive semi-definite. The notations $\mathcal{R}(\mathbf{A})$ is used to indicate the numerical range of $\mathbf{A}$. In this paper, $\mathbb{E}$ specifically denotes expectation with respect to the random test matrix.

The rest of this paper is organized as follows. In Section 2, we provide an overview of recent randomized algorithms for low-rank matrix approximation and highlight the existing bottlenecks. Section 3 presents low-rank approximation algorithms with unknown rank, accompanied by an analysis of the algorithmic complexity. In Section 4, a detailed error analysis of the proposed algorithm is presented. Section 5 includes several simulation experiments and the application of our algorithm to image reconstruction.

\section{Related works}
This section studies some existing randomization methods for low-rank matrix approximation.

In recent years, there has been significant attention focused on randomized low-rank approximation algorithms. This interest is attributed to their computational efficiency and ease of parallel implementation, allowing for utilization on advanced computing platforms.
Halko et al. \cite{hmt11} proposed a modular framework utilizing random sampling schemes to construct random algorithms for matrix decomposition. They emphasized the use of random sampling to identify subspaces that preserve most of the action of the matrix. The input matrix was compressed to this subspace, after which deterministic operations were performed on the reduced matrix to obtain the desired low-rank decomposition. Their proposed algorithm (randomized SVD) \cite{hmt11} used a random matrix to project the input matrix onto a low-dimensional subspace, effectively capturing the basic characteristics of the matrix. Subsequent calculations involved QR decomposition and SVD of the reduced-size matrix, resulting in the desired low-rank approximation decomposition. They then proposed two-sided randomized SVD (TSR-SVD) \cite{hmt11}, which ran as a single-pass method and required only one pass of the data to compute the low-rank approximation decomposition. They confirmed through extensive numerical experiments and comprehensive error analysis that the two proposed algorithms outperformed traditional competitors in terms of accuracy, speed, and robustness.

\begin{algorithm}[!ht]
\caption{Randomized SVD \cite{hmt11}}\label{alg:siam}
    \renewcommand{\algorithmicrequire}{\textbf{Input:}}
    \renewcommand{\algorithmicensure}{\textbf{Output:}}
\begin{algorithmic}[1]
\REQUIRE $\mathbf{A}\in\mathbb{C}^{m\times n}$, integer $r\leq d\leq n$.
\ENSURE A low-rank approximation: $\tilde{\mathbf{A}}=\mathbf{U}\boldsymbol{\Sigma} \mathbf{V}^{\rmH}$, where $\boldsymbol{\Sigma}\in\mathbb{R}^{d\times d}$ is a diagonal matrix, $\mathbf{U}\in\mathbb{C}^{m\times d}$ and $\mathbf{V}\in\mathbb{C}^{n\times d}$ have orthonormal columns.   
    \STATE Generate a standard Gaussian random matrix $\boldsymbol{\Omega}\in\mathbb{R}^{n\times d}$.
    \STATE Form $\mathbf{B}=\mathbf{A}\boldsymbol{\Omega}$.
    \STATE Compute QR decomposition $\mathbf{B}=\mathbf{Q}\mathbf{R}$.
    \STATE Form $\mathbf{C}=\mathbf{Q}^{\rmH}\mathbf{A}$.
    \STATE Compute SVD $\mathbf{C}=\mathbf{U}_{1}\boldsymbol{\Sigma}\mathbf{V}^{\rmH}$.
    \RETURN $\mathbf{U}=\mathbf{Q}\mathbf{U}_{1}$, $\boldsymbol{\Sigma}$ and $\mathbf{V}$.
\end{algorithmic}
\end{algorithm}

\begin{algorithm}[!ht]
\caption{Two-sided randomized SVD (TSR-SVD) \cite{hmt11}}\label{alg:tsrsvd}
    \renewcommand{\algorithmicrequire}{\textbf{Input:}}
    \renewcommand{\algorithmicensure}{\textbf{Output:}}
\begin{algorithmic}[1]
\REQUIRE $\mathbf{A}\in\mathbb{C}^{m\times n}$, integer $r\leq d\leq n$.
\ENSURE A low-rank approximation: $\tilde{\mathbf{A}}=\mathbf{U}\boldsymbol{\Sigma} \mathbf{V}^{\rmH}$, where $\boldsymbol{\Sigma}\in\mathbb{R}^{d\times d}$ is a diagonal matrix, $\mathbf{U}\in\mathbb{C}^{m\times d}$ and $\mathbf{V}\in\mathbb{C}^{n\times d}$ have orthonormal columns.   
    \STATE Generate two standard Gaussian random matrices $\boldsymbol{\Omega}_{1}\in\mathbb{R}^{n\times d}$ and $\boldsymbol{\Omega}_{2}\in\mathbb{R}^{m\times d}$.
    \STATE Form $\mathbf{B}_{1}=\mathbf{A}\boldsymbol{\Omega}_{1}$, $\mathbf{B}_{2}=\mathbf{A}^{\rmH}\boldsymbol{\Omega}_{2}$.
    \STATE Compute QR decompositions $\mathbf{B}_{1}=\mathbf{Q}_{1}\mathbf{R}_{1}$, $\mathbf{B}_{2}=\mathbf{Q}_{2}\mathbf{R}_{2}$.
    \STATE Form $\mathbf{C}=\mathbf{Q}_{1}^{\rmH}\mathbf{A}\mathbf{Q}_{2}$.
    \STATE Compute SVD $\mathbf{C}=\mathbf{U}_{1}\boldsymbol{\Sigma} \mathbf{V}_{1}^{\rmH}$.
    \RETURN $\mathbf{U}=\mathbf{Q}_{1}\mathbf{U}_{1}$, $\boldsymbol{\Sigma}$ and $\mathbf{V}=\mathbf{Q}_{2}\mathbf{V}_{1}^{H}$.
\end{algorithmic}
\end{algorithm}

Kaloorazi and Lamare \cite{kl18} proposed a compressed randomized UTV (CoR-UTV) decomposition with rank revelation, and its variant was obtained using power method techniques, as outlined in Algorithm \ref{alg:corutv}. 
\begin{definition}[UTV decomposition] \cite{kl18}
  Given a matrix $\mathbf{A}\in\mathbb{C}^{m\times n}$, then it has the UTV decomposition
  \begin{equation*}
    \mathbf{A} = \mathbf{U}\mathbf{T}\mathbf{V}^{H},
  \end{equation*}
where $\mathbf{U}\in\mathbb{C}^{m\times n}$ and $\mathbf{V}\in\mathbb{C}^{n\times n}$ have orthonormal columns, and $\mathbf{T}$ is triangular. If $\mathbf{T}$ is upper triangular, the decomposition is called URV decomposition:
\begin{equation*}
  \mathbf{A} = \mathbf{U}\left(
         \begin{array}{cc}
           \mathbf{T}_{11} & \mathbf{T}_{12} \\
           \mathbf{O} & \mathbf{T}_{22} \\
         \end{array}
       \right)
  \mathbf{V}^{H}.
\end{equation*}
If $\mathbf{T}$ is lower triangular, the decomposition is called ULV decomposition:
\begin{equation*}
  \mathbf{A} = \mathbf{U}\left(
         \begin{array}{cc}
           \mathbf{T}_{11} & \mathbf{O} \\
           \mathbf{T}_{21} & \mathbf{T}_{22} \\
         \end{array}
       \right)
  \mathbf{V}^{H}.
\end{equation*}
\end{definition}
The CoR-UTV algorithm was primarily designed for approximating low-rank input matrices through the utilization of randomized sampling schemes. CoR-UTV necessitates multiple passes over the data and operates with $O(mnr)$ floating-point operations. Furthermore, the algorithm was amenable to optimization on contemporary computing platforms for achieving maximum efficiency.

\begin{algorithm}[!ht]
\caption{Compressed randomized UTV (CoR-UTV) \cite{kl18}}\label{alg:corutv}
    \renewcommand{\algorithmicrequire}{\textbf{Input:}}
    \renewcommand{\algorithmicensure}{\textbf{Output:}}
\begin{algorithmic}[1]
\REQUIRE $\mathbf{A}\in\mathbb{C}^{m\times n}$, integer $r\leq d\leq n$.
\ENSURE A low-rank approximation: $\tilde{\mathbf{A}}=\mathbf{U}\mathbf{T} \mathbf{V}^{\rmH}$, where $\mathbf{T}\in\mathbb{R}^{d\times d}$ is an upper triangular matrix, $\mathbf{U}\in\mathbb{C}^{m\times d}$ and $\mathbf{V}\in\mathbb{C}^{n\times d}$ have orthonormal columns.   
    \STATE Generate a standard Gaussian random matrix $\boldsymbol{\Omega}\in\mathbb{R}^{n\times d}$.
    \STATE Form $\mathbf{B}_{1}=\mathbf{A}\boldsymbol{\Omega}$.
    \STATE Form $\mathbf{B}_{2}=\mathbf{A}^{\rmH}\mathbf{B}_{1}$.
    \STATE Compute QR decompositions $\mathbf{B}_{1}=\mathbf{Q}_{1}\mathbf{R}_{1}$, $\mathbf{B}_{2}=\mathbf{Q}_{2}\mathbf{R}_{2}$.
    \STATE Form $\mathbf{C}=\mathbf{Q}_{1}^{\rmH}\mathbf{A}\mathbf{Q}_{2}$.
    \STATE Compute QR decomposition $\mathbf{C}=\tilde{\mathbf{Q}}\tilde{\mathbf{R}}\tilde{\mathbf{\Pi}}$.
    \RETURN $\mathbf{U}=\mathbf{Q}_{1}\tilde{\mathbf{Q}}$, $\mathbf{T}=\tilde{\mathbf{R}}$ and $\mathbf{V}=\mathbf{Q}_{2}\tilde{\mathbf{\Pi}}$.
\end{algorithmic}
\end{algorithm}

Kaloorazi and Chen \cite{kc20} proposed the rank revealing randomized pivoted two-sided orthogonal decomposition (RP-TSOD) algorithm. This algorithm transformed the input matrix into a lower-dimensional space using a random matrix. The column-pivoted QR (CPQR) algorithm was then applied to the corresponding reduced-size matrix to construct the decomposition. This process yielded the (approximately) dominant left and right singular bases and their associated singular values of the input matrix.
Three types of bounds for RP-TSOD were provided: (i) upper bounds on the error of the low-rank approximation, (ii) error bounds for the $r$ approximate principal singular values, and (iii) bounds for the canonical angles between the approximate and exact singular subspaces.

\begin{algorithm}[!ht]
\caption{Randomized pivoted two-sided orthogonal decomposition (RP-TSOD) \cite{kc20}}\label{alg:rpstsod}
    \renewcommand{\algorithmicrequire}{\textbf{Input:}}
    \renewcommand{\algorithmicensure}{\textbf{Output:}}
\begin{algorithmic}[1]
\REQUIRE $\mathbf{A}\in\mathbb{C}^{m\times n}$, integer $r\leq d\leq n$.
\ENSURE A low-rank approximation: $\tilde{\mathbf{A}}=\mathbf{U}\mathbf{D} \mathbf{V}^{\rmH}$, where $\mathbf{D}\in\mathbb{R}^{d\times d}$ is a lower triangular matrix, $\mathbf{U}\in\mathbb{C}^{m\times d}$ and $\mathbf{V}\in\mathbb{C}^{n\times d}$ have orthonormal columns.   
    \STATE Generate a standard Gaussian random matrix $\boldsymbol{\Omega}\in\mathbb{R}^{d\times m}$.
    \STATE Form $\mathbf{B}=\boldsymbol{\Omega}\mathbf{A}$.
    \STATE Compute QR decomposition $\mathbf{B}^{\rmH}=\mathbf{Q}\mathbf{R}$.
    \STATE Form $\mathbf{C}=\mathbf{A}\mathbf{Q}$.
    \STATE Perform CPQR $\mathbf{C}=\bar{\mathbf{Q}}\bar{\mathbf{R}}\bar{\mathbf{\Pi}}^{\rmH}$.
    \STATE Perform CPQR $\bar{\mathbf{R}}^{\rmH}=\hat{\mathbf{Q}}\hat{\mathbf{R}}\hat{\mathbf{\Pi}}^{\rmH}$.
    \RETURN $\mathbf{U}=\bar{\mathbf{Q}}\hat{\mathbf{\Pi}}$, $\mathbf{D}=\hat{\mathbf{R}}^{\rmH}$ and $\mathbf{V}=\mathbf{Q}\bar{\mathbf{\Pi}}\hat{\mathbf{Q}}$.
\end{algorithmic}
\end{algorithm}

Kaloorazi and Chen \cite{kc21} introduced the projection-based partial QLP (PbP-QLP), which provided an efficient approximation of the ULV decomposition while maintaining high accuracy.
The algorithm was based on randomization. In comparison to ULV, PbP-QLP dis not use the pivoting strategy. Consequently, PbP-QLP could harness modern computer architectures more effectively, potentially surpassing competing randomized algorithms.

\begin{algorithm}[!ht]
\caption{Projection-based partial QLP (PbP-QLP) \cite{kc21}}\label{alg:pbpqlp}
    \renewcommand{\algorithmicrequire}{\textbf{Input:}}
    \renewcommand{\algorithmicensure}{\textbf{Output:}}
\begin{algorithmic}[1]
\REQUIRE $\mathbf{A}\in\mathbb{C}^{m\times n}$, integer $r\leq d\leq n$.
\ENSURE A low-rank approximation: $\tilde{\mathbf{A}}=\mathbf{Q}\mathbf{L} \mathbf{P}^{\rmH}$, where $\mathbf{L}\in\mathbb{R}^{d\times d}$ is a lower triangular matrix, $\mathbf{Q}\in\mathbb{C}^{m\times d}$ and $\mathbf{P}\in\mathbb{C}^{n\times d}$ have orthonormal columns.   
    \STATE Generate a standard Gaussian random matrix $\boldsymbol{\Omega}\in\mathbb{R}^{m\times d}$.
    \STATE Form $\mathbf{B}=\mathbf{A}^{\rmH}\boldsymbol{\Omega}$.
    \STATE Compute  an orthonormal basis for $\mathcal{R}(\mathbf{B})$: $\tilde{\mathbf{Q}}=\texttt{orth}(\mathbf{B})$.
    \STATE Form $\mathbf{C}=\mathbf{A}\tilde{\mathbf{Q}}$.
    \STATE Compute QR decomposition $\mathbf{C}=\mathbf{Q}\mathbf{R}$.
    \STATE Compute QR decomposition $\mathbf{R}^{\rmH}=\hat{\mathbf{Q}}\hat{\mathbf{R}}$.
    \RETURN $\mathbf{Q}$, $\mathbf{L}=\hat{\mathbf{R}}^{\rmH}$ and $\mathbf{P}=\tilde{\mathbf{Q}}\hat{\mathbf{Q}}$.
\end{algorithmic}
\end{algorithm}

When $\tau = 0$, all the above algorithms may be sufficiently accurate for matrices whose singular values display some decay, however in applications where the data matrix has slowly decaying singular values, it may produce singular vectors and singular values that deviate significantly from the exact ones.
Thus, variants of all the above algorithms incorporate $q$ steps of a power iteration to improve the accuracy of the algorithms in these circumstances.

Most of the existing random methods set a numerical rank $r$ in advance. However, the selection of $r$ is often obtained based on experience in practice. This is the bottleneck that this paper aims to break through. We will design an automatic basis extraction random algorithm to achieve rank $r$ approximation of low-rank matrices.


\section{Efficient orthogonal decomposition with Automatic Basis Extraction}
In this section we construct randomized algorithm for low-rank matrix $\mathbf{A}$ as in \eqref{EQLP}. Gaussian random matrices are used to construct the randomized algorithm.
We will consider a matrix $\mathbf{A}\in\mathbb{C}^{m\times n}$ with $m\geq n$, and all results are easily extended to the case of $m<n$.

\subsection{Randomized algorithms for basis extraction}
Given a matrix $\mathbf{A}$, when its rank is unknown, obviously none of the above algorithms can obtain the optimal rank $r$ approximation of matrix $\mathbf{A}$, so this section first designs a random algorithm for basic extraction that automatically finds the rank of the matrix, see Algorithm \ref{alg:qr}.

Since the size of the approximated basis of $\mathbf{U}$ is unknown in a priori, we conduct an iterative scheme to obtain the basis batch by batch. We could also calculate the basis one by one, which is less efficient in modern computer architecture. Hence, we define a blocksize hyperparameter $k$ controlling the batch size to benefit from the memory hierarchy efficiency. It is mentioned in \cite{p21} that in many environments, picking $k$ between 10 and 100 would be about right. 

For each iteration in the basis extraction, we apply the matrix to a standard Gaussian random matrix of size $n \times k$, followed by a projection matrix projecting out the bases from previous iterations. Then, we apply the QR decomposition to the matrix product result and obtain another batch of bases. We seek to build an orthonormal matrix $\mathbf{Q}$ such that
$$\|(\mathbf{I} - \mathbf{Q}\mathbf{Q}^H) \mathbf{A}\| <\varepsilon.$$
This stopping criterion can be translated into the case of judging the diagonal elements of the upper triangular matrix $\mathbf{R}$.
Hence, the dominant computational cost of the basis extraction algorithm lies in applying the matrix $\mathbf{A}$ to standard Gaussian random matrices.

\begin{algorithm}[!ht]
    \caption{Randomized algorithms for basis extraction}
    \label{alg:qr}
    \renewcommand{\algorithmicrequire}{\textbf{Input:}}
    \renewcommand{\algorithmicensure}{\textbf{Output:}}
    \begin{algorithmic}[1]
    \REQUIRE Low-rank matrix $\mathbf{A}\in\mathbb{C}^{m\times n}$, standard Gaussian matrix $\boldsymbol{\Omega}=(\boldsymbol{\Omega}_{1},\boldsymbol{\Omega}_{2},\ldots,\boldsymbol{\Omega}_{\frac{n}{k}})\in\mathbb{R}^{n\times n}$, $1\leq k\leq n$ and precision $\varepsilon$.  
    \ENSURE An approximated basis $\mathbf{Q}\in\mathbb{C}^{m\times r}$ of $\mathbf{U}$ for $\mathbf{A} = \mathbf{U}\mathbf{D} \mathbf{V}^{H}$ as in \eqref{EQLP}.   

%
%
%

\FOR {$j=1,2,\ldots,\frac{n}{k}$} \label{for}

\STATE Compute QR decomposition $(\mathbf{I}_{m}-\sum\limits_{i=1}^{j-1}\mathbf{Q}_{i}\mathbf{Q}_{i}^{\rmH})\mathbf{A}\boldsymbol{\Omega}_{j} = \mathbf{Q}_{j}\mathbf{R}_{j}$, where $\mathbf{Q}_{j} \in \mathbb{C}^{m\times k}$, $\mathbf{R}_{j} \in \mathbb{C}^{k\times k}$.

\IF{find the minimum $l$ with $|\mathbf{R}_{j}(l,l)|<\varepsilon$}

\STATE $\mathbf{Q}=[\mathbf{Q}_{1},\mathbf{Q}_{2},\ldots,\mathbf{Q}_{j}(:,1:l-1)]$.

\STATE break.

\ENDIF

\ENDFOR

        \RETURN $\mathbf{Q}$. \label{return}
    \end{algorithmic}
\end{algorithm}

Next, we explain the rationality of the stopping criterion in step 3 in Algorithm \ref{alg:qr}.
\begin{lemma}\label{lem1}\cite{hmt11}
Let $\boldsymbol{\Omega}\in \mathbb{R}^{m\times n}$ be standard Gaussian matrix. Let $\mathbf{U}=\left(
                                                                                                            \begin{array}{c}
                                                                                                              \mathbf{U}_{1} \\
                                                                                                              \mathbf{U}_{2} \\
                                                                                                            \end{array}
                                                                                                          \right)
\in\mathbb{C}^{m\times m}$ and $\mathbf{V}=(\mathbf{V}_{1},\mathbf{V}_{2})\in\mathbb{C}^{n\times n}$ be orthonormal matrices. Let $\mathbf{U}\boldsymbol{\Omega}=\left(
                                                                                                            \begin{array}{c}
                                                                                                              \mathbf{U}_{1}\boldsymbol{\Omega} \\
                                                                                                              \mathbf{U}_{2}\boldsymbol{\Omega} \\
                                                                                                            \end{array}
                                                                                                          \right) 
=\left(
                                                                                                            \begin{array}{c}
                                                                                                              \boldsymbol{\Omega}_{1} \\
                                                                                                              \boldsymbol{\Omega}_{2} \\
                                                                                                            \end{array}
                                                                                                          \right)$ and
$\boldsymbol{\Omega} \mathbf{V}=\boldsymbol{\Omega}(\mathbf{V}_{1},\mathbf{V}_{2})=(\boldsymbol{\Omega} \mathbf{V}_{1},\boldsymbol{\Omega} \mathbf{V}_{2})=(\mathbf{\Lambda}_{1},\mathbf{\Lambda}_{2})$. Then we have the following
conclusions.
\begin{enumerate}
  \item If $m \leq n$, then $\boldsymbol{\Omega}$ is a full row rank with probability one. If $m \geq n$,
then $\boldsymbol{\Omega}$ is a full column rank with probability one. Then $rank(\boldsymbol{\Omega}) =
\min\{m, n\}$.
  \item $\mathbf{U}\boldsymbol{\Omega}$, $\boldsymbol{\Omega} \mathbf{V}$, $\boldsymbol{\Omega}_{1}$, $\boldsymbol{\Omega}_{2}$, $\mathbf{\Lambda}_{1}$ and $\mathbf{\Lambda}_{2}$ are also standard Gaussian matrices.
\end{enumerate}
\end{lemma}

\begin{lemma}\label{lem2}
Let $\mathbf{A} \in \mathbb{C}^{m\times n} (m\geq n)$ with $rank(\mathbf{A})=r$ , $\boldsymbol{\Omega}=(\boldsymbol{\Omega}_1,\boldsymbol{\Omega}_2)\in\mathbb{R}^{n\times n}$ be standard Gaussian matrix, where $\boldsymbol{\Omega}_1 \in \mathbb{R}^{n\times l}$, $\boldsymbol{\Omega}_2 \in \mathbb{R}^{n\times (n-l)}$. Then, $rank(\mathbf{A}\boldsymbol{\Omega}_1)=\min\{r,l\}$.
\end{lemma}

\begin{IEEEproof}
Since the SVD of $\mathbf{A}$ is $\mathbf{A}=\mathbf{U}\boldsymbol{\Sigma} \mathbf{V}^{\rmH}$, where $\mathbf{U} \in \mathbb{C}^{m \times m}$ and $\mathbf{V} \in \mathbb{C}
^{n\times n}$ are unitary matrices, $\boldsymbol{\Sigma} \in \mathbb{R}^{m \times n}$ is a diagonal matrix. For $rank(\mathbf{A})=r$, we have
\begin{eqnarray*}
\mathbf{A}\boldsymbol{\Omega}_1&=&\mathbf{U}\boldsymbol{\Sigma} \mathbf{V}^{\rmH}\boldsymbol{\Omega}_1\\
&=&\left(
                              \begin{array}{cc}
                                \mathbf{U}_1 & \mathbf{U}_2 \\
                              \end{array}
                            \right)
                            \left(
                               \begin{array}{cc}
                                 \boldsymbol{\Sigma}_1 & \mathbf{O} \\
                                 \mathbf{O}  & \mathbf{O}\\
                               \end{array}
                             \right)
                             \left(
                               \begin{array}{c}
                                 \mathbf{V}_1^{\rmH} \\
                                 \mathbf{V}_2^{\rmH} \\
                               \end{array}
                             \right)
                             \boldsymbol{\Omega}_1\\
                             &=&\mathbf{U}_1\boldsymbol{\Sigma}_1\mathbf{V}_1^{\rmH}\boldsymbol{\Omega}_1,
\end{eqnarray*}
where $\mathbf{U}_1 \in \mathbb{U}^{m \times r}$ and $\mathbf{V}_1\in \mathbb{U}^{n \times r}$ are column orthonormal matrices, $\boldsymbol{\Sigma}_1\in \mathbb{R}^{r \times r}$ is a diagonal matrix. By Lemma \ref{lem1}, $\mathbf{V}_1^{H}\boldsymbol{\Omega}_1$ is an $r \times l$ standard Gaussian matrix, $rank(\mathbf{V}_1^{\rmH}\boldsymbol{\Omega}_1)=\min\{r,l\}$. And $\mathbf{U}_1\boldsymbol{\Sigma}_1$ is a full column rank matrix. Therefore $rank(\mathbf{A}\boldsymbol{\Omega}_1)=rank(\mathbf{U}_1\boldsymbol{\Sigma}_1\mathbf{V}_1^{\rmH}\boldsymbol{\Omega}_1)=rank(\mathbf{V}_1^{\rmH}\boldsymbol{\Omega}_1)=\min\{r,l\}$.
This completes the proof.
\end{IEEEproof}

\begin{theorem}
Let $\mathbf{A} \in \mathbb{C}^{m\times n} (m\geq n)$ with $rank(\mathbf{A})=r$, $\boldsymbol{\Omega}=(\boldsymbol{\Omega}_1,\boldsymbol{\Omega}_2)\in\mathbb{R}^{n\times n}$ be standard Gaussian matrix, where $\boldsymbol{\Omega}_1 \in \mathbb{R}^{n\times l}, \boldsymbol{\Omega}_2 \in \mathbb{R}^{n\times (n-l)}$ and $l>r$. $\mathbf{A}\boldsymbol{\Omega}_1$ has the QR decomposition $\mathbf{A}\boldsymbol{\Omega}_1=\mathbf{Q}\mathbf{R}$. Then $\mathbf{R}(i,i)=0, i=r+1,...,l$.
\end{theorem}

\begin{IEEEproof}
Observe that $\mathbf{A}\boldsymbol{\Omega}$ has the QR decomposition $\mathbf{A}\boldsymbol{\Omega}=\mathbf{P}\mathbf{T}$, where $\mathbf{P}=(\mathbf{Q}_1,\mathbf{Q}_2)$, $\mathbf{T}=\left(
        \begin{array}{cc}
          \mathbf{R}_{11} & \mathbf{R}_{12} \\
          \mathbf{O} & \mathbf{R}_{22} \\
        \end{array}
      \right)
$, $\mathbf{Q}_1 \in \mathbb{C}^{m\times l}$, $\mathbf{Q}_2 \in \mathbb{C}^{m\times (m-l)}$, $\mathbf{R}_{11} \in \mathbb{C}^{l\times l}$, $\mathbf{R}_{12} \in \mathbb{C}^{l\times (n-l)}$ and $\mathbf{R}_{22} \in \mathbb{C}^{(m-l)\times (n-l)}$. Thus $\mathbf{A}\boldsymbol{\Omega}_1=\mathbf{Q}_1\mathbf{R}_{11}$. When $l\leq r$, by Lemma \ref{lem2}, $rank(\mathbf{A}\boldsymbol{\Omega}_1)=l$, then $rank(\mathbf{R}_{11})=l$. It means that $|\mathbf{T}(j,j)|>0$, $j=1,2,\ldots,r$. When $l>r$, by Lemma \ref{lem2}, $rank(\mathbf{A}\boldsymbol{\Omega}_1)=r$. Then $rank(\mathbf{R}_{11})=r$. Since $|\mathbf{T}(j,j)|>0$, $j=1,2,\ldots,r$, $\mathbf{R}_{11} \in \mathbb{C}^{l\times l}$, we have $\mathbf{R}_{11}(i,i)=0, i=r+1,...,l$.
This completes the proof.
\end{IEEEproof}

\begin{theorem}
Let $\mathbf{A} \in \mathbb{C}^{m\times n} (m \ge n)$ with $rank(\mathbf{A})=r$, $\boldsymbol{\Omega}\in \mathbb{R}^{n\times n}$ be standard Gaussian matrix, partition $\boldsymbol{\Omega}=(\boldsymbol{\Omega}_{1},\cdots,\boldsymbol{\Omega}_{s})$ into blocks containing, respectively, where $\boldsymbol{\Omega}_{j} (j=1,2,...,s)$ have the same columns. $\mathbf{A}\boldsymbol{\Omega}$ has the QR decomposition $\mathbf{A}\boldsymbol{\Omega}=\mathbf{A}(\boldsymbol{\Omega}_{1},\cdots,\boldsymbol{\Omega}_{s})=(\mathbf{Q}_1,\cdots,\mathbf{Q}_s)\left(
                                                   \begin{array}{ccc}
                                                     \mathbf{R}_{11} & \cdots & \mathbf{R}_{1s} \\
                                                      & \ddots & \vdots \\
                                                      &  & \mathbf{R}_{ss} \\
                                                   \end{array}
                                                 \right)$.
Let $\mathbf{Y}_j=\mathbf{A}\boldsymbol{\Omega}_{j}-\sum_{i=1}^{j-1}\mathbf{Q}_{i}\mathbf{Q}_{i}^{\rmH}\mathbf{A}\boldsymbol{\Omega}_{j}$, then $\mathbf{Y}_j$ has the QR decomposition $\mathbf{Y}_j=\mathbf{P}_j\mathbf{T}_j$, where $\mathbf{P}_{j}=\mathbf{Q}_j$, $\mathbf{T}_j=\mathbf{R}_{jj}$, $j=1,2,\ldots,s$.
\end{theorem}
\begin{IEEEproof}
When $s=1$, $\mathbf{Y}_1=\mathbf{A}\boldsymbol{\Omega}_{1}$, $\mathbf{Y}_1$ has the QR decomposition $\mathbf{Y}_1=\mathbf{P}_{1}\mathbf{T}_{1}$. Due to $\mathbf{Y}_1$ is a full column rank matrix, the QR decomposition of $\mathbf{Y}_1$ is unique. Observe that $\mathbf{P}_{1}=\mathbf{Q}_1$, $\mathbf{T}_1=\mathbf{R}_{11}$.
When $s=2$, $\mathbf{Y}_2=\mathbf{A}\boldsymbol{\Omega}_{2}-\mathbf{Q}_{1}\mathbf{Q}_{1}^{\rmH}\mathbf{A}\boldsymbol{\Omega}_{2}$, $\mathbf{Y}_2$ has the QR decomposition $\mathbf{Y}_2=\mathbf{P}_{2}\mathbf{T}_{2}$. For $\mathbf{A}(\boldsymbol{\Omega}_{1},\boldsymbol{\Omega}_{2})=(\mathbf{Q}_{1},\mathbf{Q}_{2})
\left(
  \begin{array}{cc}
    \mathbf{R}_{11} & \mathbf{R}_{12} \\
    \mathbf{O} & \mathbf{R}_{22} \\
  \end{array}
\right)$.
It follows that
\begin{equation*}
\mathbf{A}\boldsymbol{\Omega}_{1}=\mathbf{Q}_{1}\mathbf{R}_{11}, \mathbf{A}\boldsymbol{\Omega}_{2}=\mathbf{Q}_{1}\mathbf{R}_{12}+\mathbf{Q}_{2}\mathbf{R}_{22}
\end{equation*}
and
\begin{equation*}
\begin{aligned}
\left(                                                                                        \begin{array}{cc}
                                                                                          \mathbf{R}_{11} & \mathbf{R}_{12} \\
                                                                                          \mathbf{O} & \mathbf{R}_{22} \\
                                                                                        \end{array}
                                                                                      \right)
                                                                                              &=\left(
                                                                                                \begin{array}{c}
                                                                                                  \mathbf{Q}_1^H \\
                                                                                                  \mathbf{Q}_2^H \\
                                                                                                \end{array}
                                                                                              \right)\left(
                                                                                                                    \begin{array}{cc}
                                                                                                                      \mathbf{Q}_{1} & \mathbf{Q}_{2} \\
                                                                                                                    \end{array}
                                                                                                                  \right)
\left(
  \begin{array}{cc}
    \mathbf{R}_{11} & \mathbf{R}_{12} \\
    \mathbf{O} & \mathbf{R}_{22} \\
  \end{array}
\right)\\&=\left(
                                                                                                \begin{array}{c}
                                                                                                  \mathbf{Q}_1^H \\
                                                                                                  \mathbf{Q}_2^H \\
                                                                                                \end{array}
                                                                                              \right)
                                                                                      \left(
                                                                                                \begin{array}{cc}
                                                                                                  \mathbf{A}\boldsymbol{\Omega}_1 & \mathbf{A}\boldsymbol{\Omega}_2 \\
                                                                                                \end{array}
                                                                                              \right) \\
&=\left(
          \begin{array}{cc}
            \mathbf{Q}_1^H\mathbf{A}\boldsymbol{\Omega}_1 & \mathbf{Q}_1^H\mathbf{A}\boldsymbol{\Omega}_2 \\
            \mathbf{Q}_2^H\mathbf{A}\boldsymbol{\Omega}_1 & \mathbf{Q}_2^H\mathbf{A}\boldsymbol{\Omega}_2 \\
          \end{array}
        \right).
        \end{aligned}
\end{equation*}
Observe that $\mathbf{Q}_{2}\mathbf{R}_{22}=\mathbf{A}\boldsymbol{\Omega}_{2}-\mathbf{Q}_{1}\mathbf{R}_{12}=\mathbf{A}\boldsymbol{\Omega}_{2}-\mathbf{Q}_{1}\mathbf{Q}_{1}^{\rmH}\mathbf{A}\boldsymbol{\Omega}_{2}$. Due to $\mathbf{Y}_2$ is a full column rank matrix, the QR decomposition of $\mathbf{Y}_2$ is unique. So $\mathbf{P}_{2}=\mathbf{Q}_2$, $\mathbf{T}_2=\mathbf{R}_{22}$. The same is true when $s>2$.
This completes the proof.
\end{IEEEproof}

\subsection{Efficient orthogonal decomposition with automatic basis extraction}

Next we describe an efficient orthogonal decomposition with automatic basis extraction for computing matrix of unknown rank.

Given the matrix $\mathbf{A}$, the procedure to compute the basic form of efficient orthogonal decomposition with unknown rank is as follows: we compute an approximated basis $\mathbf{Q}_{0}$ of $\mathbf{U}$ for $\mathbf{A}=\mathbf{U}\mathbf{D} \mathbf{V}^{\rmH}$ by Algorithm \ref{alg:qr}.
We then form the matrix:
\begin{equation}\label{CQA}
\mathbf{C}=\mathbf{Q}_{0}^{\rmH}\mathbf{A}.
\end{equation}
The matrix $\mathbf{C} \in \mathbb{C}^{r\times n}$ is formed by linear combinations of $\mathbf{A}$'s rows by means of $\mathbf{Q}$. Then we compute the QR decomposition of $\mathbf{C}^{H}$:
\begin{equation}\label{CQR}
\mathbf{C}^{H}=\mathbf{Q}\mathbf{R},
\end{equation}
where $\mathbf{Q} \in \mathbb{C}^{n\times r}$ is column-orthonormal, $\mathbf{R}\in \mathbb{C}^{r\times r}$ is an upper triangular matrix. We call the diagonals of $\mathbf{R}$, R-values. Next, we compute QR decomposition on $\mathbf{R}^{H}$:
\begin{equation}\label{RQR}
\mathbf{R}^{H}=\hat{\mathbf{Q}}\hat{\mathbf{R}},
\end{equation}
where $\hat{\mathbf{Q}} \in \mathbb{C}^{r\times r}$ is orthonormal matrix, $\hat{\mathbf{R}}\in \mathbb{C}^{r\times r}$ is an upper triangular matrix.
Lastly, we construct the low-rank approximation of $\mathbf{A}$:
\begin{eqnarray}\label{AUDV}
\tilde{\mathbf{A}}=\mathbf{Q}_{0}(\mathbf{Q}_{0}^{H}\mathbf{A})
=\mathbf{Q}_{0}\mathbf{R}^{H}\mathbf{Q}^{H}\nonumber\\
=\mathbf{Q}_{0}\mathbf{\hat{Q}}\mathbf{\hat{R}}\mathbf{Q}^{H} := \mathbf{U}\mathbf{D} \mathbf{V}^H,
\end{eqnarray}
where $\mathbf{U}=\mathbf{Q}_{0}\mathbf{\hat{Q}}$, $\mathbf{D}=\mathbf{\hat{R}}$ and $\mathbf{V}=\mathbf{Q}$.

The overall algorithm for computing low-rank approximation of matrix of unknown rank is summarized in Algorithm \ref{alg:alg1}, where the basis extraction algorithm (Algorithm \ref{alg:qr}) is denoted as ``BasisExt''.

\begin{algorithm}[!ht]
    \caption{Efficient orthogonal decomposition with automatic basis extraction}
    \label{alg:alg1}
    \renewcommand{\algorithmicrequire}{\textbf{Input:}}
    \renewcommand{\algorithmicensure}{\textbf{Output:}}
    \begin{algorithmic}[1]
        \REQUIRE $\mathbf{A}\in\mathbb{C}^{m\times n}$ and precision $\varepsilon$.  
        \ENSURE Rank $r$ approximation: $\tilde{\mathbf{A}}=\mathbf{U}\mathbf{D} \mathbf{V}^{\rmH}$, where $\mathbf{D}\in\mathbb{R}^{r\times r}$ is an upper triangular matrix, $\mathbf{U}\in\mathbb{C}^{m\times r}$ and $\mathbf{V}\in\mathbb{C}^{n\times r}$ have orthonormal columns.  
%
        \STATE Generate a standard Gaussian matrix $\boldsymbol{\Omega}\in\mathbb{R}^{n\times n}$.
        \STATE Compute $\mathbf{Q}_{0}=BasisExt(\mathbf{A},\boldsymbol{\Omega}, \varepsilon)$ by Algorithm \ref{alg:qr}.
        \STATE Form $\mathbf{C}=\mathbf{Q}_{0}^{\rmH}\mathbf{A}$.
        \STATE Compute QR decomposition $\mathbf{C}^{\rmH}=\mathbf{Q}\mathbf{R}$.
        \STATE Compute QR decomposition $\mathbf{R}^{H}=\mathbf{\hat{Q}}\mathbf{\hat{R}}$.

        \RETURN $\mathbf{U}=\mathbf{Q}_{0}\mathbf{\hat{Q}}$, $\mathbf{D}=\mathbf{\hat{R}}$ and $\mathbf{V}=\mathbf{Q}$.
    \end{algorithmic}
\end{algorithm}

\subsection{Subspace iteration - EOD-ABE}
When the input matrix $\mathbf{A}$ has flat singular values or is very large, the influence of these small singular values relative to the main singular vectors can also be reduced by exponentiating the input matrix, as detailed in Algorithm \ref{alg:subiter}.
In the first step of Algorithm \ref{alg:subiter}, the rank of matrix $\mathbf{A}$ can be obtained, which is the number of columns of $\mathbf{Q}$. Therefore, in steps 4 and 6 of the algorithm, only QR decomposition is required, and there is no need to use Algorithm \ref{alg:qr} to find the basis.

\begin{algorithm}[H]
    \caption{Orthorgonalization with QR}
    \label{alg:subiter}
    \renewcommand{\algorithmicrequire}{\textbf{Input:}}
    \renewcommand{\algorithmicensure}{\textbf{Output:}}
    \begin{algorithmic}[1]
        \REQUIRE $\mathbf{A}\in\mathbb{C}^{m\times n}$, a standard Gaussian matrix $\boldsymbol{\Omega}\in\mathbb{R}^{n\times n}$, precision $\varepsilon$ and a power iteration $\tau$.  
        \ENSURE An approximated basis $\mathbf{Q}_{\tau}$ for $\mathbf{Y}_{\tau}=(\mathbf{A}\mathbf{A}^{H})^{\tau}\mathbf{A}\boldsymbol{\Omega}$.  
        \STATE Compute $\mathbf{Q}_{0}=BasisExt(\mathbf{A},\boldsymbol{\Omega},\varepsilon)$ by Algorithm \ref{alg:qr}.
        \FOR {$j=1,2,\ldots,\tau$}
            \STATE Form $\tilde{\mathbf{Y}}_{j}=\mathbf{A}^{\rmH}\mathbf{Q}_{j-1}$.
            \STATE Compute QR decomposition $\tilde{\mathbf{Y}}_{j}=\tilde{\mathbf{Q}}_{j}\tilde{\mathbf{R}}_{j}$.
            \STATE Form $\mathbf{Y}_{j}=\mathbf{A}\tilde{\mathbf{Q}}_{j}$.
            \STATE Compute QR decomposition $\mathbf{Y}_{j}=\mathbf{Q}_{j}\mathbf{R}_{j}$.
        \ENDFOR
    \end{algorithmic}
\end{algorithm}

The efficient orthogonal decomposition for basis extraction via subspace iteration, see Algorithm \ref{alg:A} for details.

\begin{algorithm}[!ht]
    \caption{Efficient orthogonal decomposition with automatic basis extraction via subspace iteration}
    \label{alg:A}
    \renewcommand{\algorithmicrequire}{\textbf{Input:}}
    \renewcommand{\algorithmicensure}{\textbf{Output:}}
    \begin{algorithmic}[1]
        \REQUIRE $\mathbf{A}\in\mathbb{C}^{m\times n}$, precision $\varepsilon$ and a power iteration $\tau$.  
        \ENSURE Rank $r$ approximation: $\tilde{\mathbf{A}}=\mathbf{U}\mathbf{D} \mathbf{V}^{\rmH}$, where $\mathbf{D}\in\mathbb{R}^{r\times r}$ is an upper triangular matrix, $\mathbf{U}\in\mathbb{C}^{m\times r}$ and $\mathbf{V}\in\mathbb{C}^{n\times r}$ have orthonormal columns.  
        \STATE Generate a standard Gaussian matrix $\boldsymbol{\Omega}\in\mathbb{R}^{n\times n}$.
        \STATE Compute $\mathbf{Y}_{\tau}=(\mathbf{A}\mathbf{A}^{H})^{\tau}\mathbf{A}\boldsymbol{\Omega}$ and compute an approximated basis $\mathbf{Q}_{\tau}$ for $\mathbf{Y}_{\tau}$ by Algorithm \ref{alg:subiter}.
        
        \STATE Form $\mathbf{C}=\mathbf{Q}_{\tau}^{\rmH}\mathbf{A}$.
        \STATE Compute QR decomposition $\mathbf{C}^{\rmH}=\mathbf{Q}\mathbf{R}$.
        \STATE Compute QR decomposition $\mathbf{R}^{H}=\mathbf{\hat{Q}}\mathbf{\hat{R}}$.

        \RETURN $\mathbf{U}=\mathbf{Q}_{\tau}\mathbf{\hat{Q}}$, $\mathbf{D}=\mathbf{\hat{R}}$ and $\mathbf{V}=\mathbf{Q}$.
    \end{algorithmic}
\end{algorithm}

\subsection{Computational cost}

\textbf{Arithmetic Cost.} To compute an approximation of matrix $\mathbf{A}$, EOD-ABE requires the following arithmetic operations:

\begin{itemize}
  \item Computing $\mathbf{Q}_{0}$ by Algorithm \ref{alg:qr} costs $mnr+mr^2+mrk$.
  \item Forming $\mathbf{C}$ \eqref{CQA} costs $mnr$.
  \item Computing $\mathbf{Q}$ and $\mathbf{R}$ in \eqref{CQR} costs $2nr^2$.
  \item Computing $\hat{\mathbf{Q}}$ and $\hat{\mathbf{R}}$ in \eqref{RQR} costs $2r^3$.
  \item Computing $\mathbf{U}$, $\mathbf{D}$ and $\mathbf{V}$ in \eqref{AUDV} costs $mr^2$.
\end{itemize}

The total computational complexity of the above steps is $2mnr+2mr^2+2nr^2+2r^3+mrk$.

After subspace iteration, the complexity of step 2 of Algorithm \ref{alg:A} increases by
$$\tau(2mnr+2mr^2+2nr^2).$$
Therefore, the total complexity of Algorithm \ref{alg:A} is 
$$(\tau+1)(2mnr+2mr^2+2nr^2)+2r^3+mrk.$$

\textbf{Communication Cost.} Communication costs are determined by data movement between processors operating in parallel and across various levels of the memory hierarchy.
On advanced computing devices, communication costs predominantly govern the factoring of any external storage matrix.
Hence, executing any factoring algorithm with minimal communication costs is highly desirable.
Algorithm \ref{alg:qr} computes several easily parallelizable matrix-matrix multiplications.
Additionally, during QR decomposition, calculations occur on small $m\times k$ matrices each time.
Thus, Algorithm \ref{alg:A} can be efficiently executed on high-performance computing devices.

\section{Error analysis of EOD-ABE}

This section provides a detailed account of theoretical analysis for efficient orthogonal decomposition with automatic basis extraction of low-rank matrix. Here, we develop bounds for the error of the low-rank approximation and the accuracy of $r$ principal singular values.

\begin{lemma}\label{lem:pm}\cite{hmt11}
For a given matrix $\mathbf{A}$, we write $\mathbf{P}_{A}$ for the unique orthogonal
projector with $\mathrm{range}(\mathbf{P}_{A})=\mathrm{range}(\mathbf{A})$. When $\mathbf{A}$ has
full column rank, we can express this projector explicitly
\begin{eqnarray*}
\mathbf{P}_{A} = \mathbf{A}(\mathbf{A}^{\rmH} \mathbf{A})^{-1}\mathbf{A}^{\rmH}.
\end{eqnarray*}
For a unitary matrix $\mathbf{Q}$, then
\begin{eqnarray*}
\mathbf{Q}^{\rmH}\mathbf{P}_{A}\mathbf{Q} = \mathbf{P}_{Q^{\rmH}A}.
\end{eqnarray*}
\end{lemma}

An orthogonal projector is an Hermitian matrix $\mathbf{P}$ that satisfies the polynomial $\mathbf{P}^{2}=\mathbf{P}$.

\begin{lemma}\cite{hmt11}\label{lem:p}
Suppose $\mathcal{R}(\mathbf{N}) \subset \mathcal{R}(\mathbf{M})$. Then, for each matrix $\mathbf{A}$, it holds that $\|\mathbf{P}_{N}\mathbf{A}\| \leq \|\mathbf{P}_{M}\mathbf{A}\|$ and that $\|(\mathbf{I} - \mathbf{P}_{M})\mathbf{A}\| \leq \|(\mathbf{I} - \mathbf{P}_{N})\mathbf{A}\|$.
\end{lemma}

\begin{lemma}\label{lem:sgt}\cite{hmt11}
Fix matrices $\mathbf{S}$, $\mathbf{T}$, and draw a standard Gaussian matrix $\boldsymbol{\Omega}$. Then
\begin{equation*}
  \mathbb{E}\|\mathbf{S}\boldsymbol{\Omega} \mathbf{T}\|_{\rmF}^{2} = \|\mathbf{S}\|_{\rmF}^{2}\|\mathbf{T}\|_{\rmF}^{2}.
\end{equation*}
\end{lemma}

\begin{lemma}\label{lem:ginv}\cite{hmt11}
Draw a $r \times d$ standard Gaussian matrix $\boldsymbol{\Omega}$ with $r \geq 2$ and $d-r \geq 2$. Then
\begin{equation*}
  \mathbb{E}\|\boldsymbol{\Omega}^{\dagger}\|_{\rmF}^{2} = \frac{r}{d-r-1}.
\end{equation*}
\end{lemma}



\begin{theorem}\label{siam}
Let $\mathbf{A}$ be an $m\times n (m\geq n)$ matrix with singular values $\sigma_{1}\geq\sigma_{2}\geq\cdots\geq\sigma_{n}$, $\mathrm{rank}(\mathbf{A})=r$. Let $\boldsymbol{\Omega}\in\mathbb{R}^{n\times d}$ be standard Gaussian matrix. The SVD of $\mathbf{A}$ is
\begin{eqnarray}
  \mathbf{A} &=& \mathbf{U}_{m\times n}\boldsymbol{\Sigma}_{n\times n} \mathbf{V}_{n\times n}^{\rmH} \nonumber\\
   &=& \left(
          \begin{array}{cc}
            \mathbf{U}_{r} & \mathbf{U}_{0} \\
          \end{array}
        \right)
        \left(
          \begin{array}{cc}
            \boldsymbol{\Sigma}_{r} & \mathbf{O} \\
            \mathbf{O} & \boldsymbol{\Sigma}_{0} \\
          \end{array}
        \right)
        \left(
          \begin{array}{c}
            \mathbf{V}_{r}^{\rmH} \\
            \mathbf{V}_{0}^{\rmH} \\
          \end{array}
        \right),
\end{eqnarray}
where $\mathbf{U}_{r}\in\mathbb{C}^{m\times r}$, $\mathbf{U}_{0}\in\mathbb{C}^{m\times(n-r)}$ have orthonormal columns, $\boldsymbol{\Sigma}_{r}\in\mathbb{R}^{r\times r}$ and $\boldsymbol{\Sigma}_{0}\in\mathbb{R}^{(n-r)\times(n-r)}$ are diagonal matrices, and $\mathbf{V}_{r}\in\mathbb{C}^{n\times r}$ and $\mathbf{V}_{0}\in\mathbb{C}^{n\times(n-r)}$ have orthonormal columns. 
Let $\tilde{\mathbf{A}}$ be a low-rank approximation to $\mathbf{A}$ computed through Algorithm \ref{alg:A}. $r\geq 2$ and $d-r\geq2$. Then
\begin{equation*}
  \|\mathbf{A}-\tilde{\mathbf{A}}\|_{\rmF}^{2} \leq \alpha^{2\tau}\|\boldsymbol{\Sigma}_{0}\mathbf{\Lambda}_{2}\mathbf{\Lambda}_{1}^{\dagger}\|_{\rmF}^{2}+\|\boldsymbol{\Sigma}_{0}\|_{\rmF}^{2}
\end{equation*}
and
\begin{equation*}
    \mathbb{E}\|\mathbf{A}-\tilde{\mathbf{A}}\|_{\rmF} \leq \left(1+\frac{r\alpha^{2\tau}}{d-r-1}\right)^{\frac{1}{2}}\left(\sum_{j>r}\sigma_{j}^{2}\right)^{\frac{1}{2}},
\end{equation*}
where $\alpha=\frac{\sigma_{r+1}}{\sigma_{r}}$, $\mathbf{\Lambda}_{1}=\mathbf{V}_{r}^{\rmH}\boldsymbol{\Omega}$ and $\mathbf{\Lambda}_{2}=\mathbf{V}_{0}^{\rmH}\boldsymbol{\Omega}$.
\end{theorem}

\begin{IEEEproof}
By Lemma \ref{lem:pm}, we write
\begin{equation}\label{th3.3:1}
  \|\mathbf{A}-\tilde{\mathbf{A}}\|_{\rmF}=\|\mathbf{A}-\mathbf{Q}\mathbf{Q}^{\rmH}\mathbf{A}\|_{\rmF}=\|(\mathbf{I}-\mathbf{P}_{Q})\mathbf{A}\|_{\rmF}.
\end{equation}
We observe that $\mathbf{B}$ is represented as:
\begin{equation*}
  \mathbf{B} = (\mathbf{A}\mathbf{A}^{\rmH})^{\tau}\mathbf{A}\boldsymbol{\Omega} = \mathbf{U}\left(
                   \begin{array}{cc}
                     \boldsymbol{\Sigma}_{r}^{2\tau+1} & \mathbf{O} \\
                     \mathbf{O} & \boldsymbol{\Sigma}_{0}^{2\tau+1} \\
                   \end{array}
                 \right)\mathbf{V}^{\rmH}\boldsymbol{\Omega}.
\end{equation*}
Let $\mathbf{V}^{\rmH}\boldsymbol{\Omega} = \left(
                     \begin{array}{c}
                       \mathbf{V}_{r}^{\rmH}\boldsymbol{\Omega} \\
                       \mathbf{V}_{0}^{\rmH}\boldsymbol{\Omega} \\
                     \end{array}
                   \right) :=
                   \left(
                     \begin{array}{c}
                       \mathbf{\Lambda}_{1} \\
                       \mathbf{\Lambda}_{2} \\
                     \end{array}
                   \right)
$.
Now, we form $\bar{\mathbf{B}}$ as:
\begin{equation*}
  \bar{\mathbf{B}} = \mathbf{U}^{\rmH}\mathbf{B} = \left(
                       \begin{array}{c}
                         \boldsymbol{\Sigma}_{r}^{2\tau+1}\mathbf{\Lambda}_{1} \\
                         \boldsymbol{\Sigma}_{0}^{2\tau+1}\mathbf{\Lambda}_{2} \\
                       \end{array}
                     \right).
\end{equation*}
We therefore have
\begin{equation*}
  \mathcal{R}(\bar{\mathbf{B}})=\mathcal{R}(\mathbf{U}^{\rmH}\mathbf{B})=\mathcal{R}(\mathbf{U}^{\rmH}\mathbf{Q}).
\end{equation*}
We define a matrix $\mathbf{X}$ as follows:
\begin{equation*}
  \mathbf{X} = \mathbf{\Lambda}_{1}^{\dagger}\boldsymbol{\Sigma}_{r}^{-(2\tau+1)},
\end{equation*}
and, form another matrix $\tilde{\mathbf{B}}$ by shrinking the subspace of $\bar{\mathbf{B}}$ through $\mathbf{X}$:
\begin{equation*}
  \tilde{\mathbf{B}} = \bar{\mathbf{B}}\mathbf{X} = \left(
                           \begin{array}{c}
                             \mathbf{I} \\
                             \mathbf{S} \\
                           \end{array}
                         \right),
\end{equation*}
where $\mathbf{S} = \boldsymbol{\Sigma}_{0}^{2\tau+1}\mathbf{\Lambda}_{2}\mathbf{\Lambda}_{1}^{\dagger}\boldsymbol{\Sigma}_{r}^{-(2\tau+1)}$.
By Lemma \ref{lem:p}, it follows that
\begin{equation}\label{range}
  \mathcal{R}(\tilde{\mathbf{B}}) \subset \mathcal{R}(\bar{\mathbf{B}}) = \mathcal{R}(\mathbf{U}^{\rmH}\mathbf{Q}),
\end{equation}
\begin{equation*}
  \|\mathbf{I}-\mathbf{P}_{U^{\rmH}Q}\|_{\rmF}\leq\|\mathbf{I}-\mathbf{P}_{\tilde{B}}\|_{\rmF}.
\end{equation*}
Then
\begin{eqnarray}\label{po4}
&&\mathbf{I}-\mathbf{P}_{\tilde{B}}\nonumber\\
&&=\left(\begin{array}{cc}
\mathbf{I}-(\mathbf{I}+\mathbf{S}^{\rmH}\mathbf{S})^{-1}&-(\mathbf{I}+\mathbf{S}^{\rmH}\mathbf{S})^{-1}\mathbf{S}^{\rmH}\\
-\mathbf{S}(\mathbf{I}+\mathbf{S}^{\rmH}\mathbf{S})^{-1}&\mathbf{I}-\mathbf{S}(\mathbf{I}+\mathbf{S}^{\rmH}\mathbf{S})^{-1}\mathbf{S}^{\rmH}
\end{array}\right).
\end{eqnarray}
Since
\begin{eqnarray*}
&&[\mathbf{I}-(\mathbf{I}+\mathbf{S}^{\rmH}\mathbf{S})^{-1}](\mathbf{S}^{\rmH}\mathbf{S})(\mathbf{I}+\mathbf{S}^{\rmH}\mathbf{S})\\
&=&\mathbf{I}+\mathbf{S}^{\rmH}\mathbf{S}-(\mathbf{I}+\mathbf{S}^{\rmH}\mathbf{S})^{-1}(\mathbf{S}^{\rmH}\mathbf{S})(\mathbf{I}+\mathbf{S}^{\rmH}\mathbf{S})\\
&=&\mathbf{I}+\mathbf{S}^{\rmH}\mathbf{S}-(\mathbf{I}+\mathbf{S}^{\rmH}\mathbf{S})^{-1}(\mathbf{I}+\mathbf{S}^{\rmH}\mathbf{S})\mathbf{S}^{\rmH}\mathbf{S}\\
&=&\mathbf{I}+\mathbf{S}^{\rmH}\mathbf{S}-\mathbf{S}^{\rmH}\mathbf{S}=\mathbf{I},
\end{eqnarray*}
then
\begin{eqnarray*}
\mathbf{I}-(\mathbf{I}+\mathbf{S}^{\rmH}\mathbf{S})^{-1}&=&(\mathbf{I}+\mathbf{S}^{\rmH}\mathbf{S})^{-1}(\mathbf{S}^{\rmH}\mathbf{S}) \\
&=&(\mathbf{S}^{\rmH}\mathbf{S})^{1/2}(\mathbf{I}+\mathbf{S}^{\rmH}\mathbf{S})^{-1}(\mathbf{S}^{\rmH}\mathbf{S})^{1/2} \\
&\preceq& \mathbf{S}^{\rmH}\mathbf{S}.
\end{eqnarray*}
It follows from (\ref{po4}) that
\begin{eqnarray}\label{eqn6}
\mathbf{I}-\mathbf{P}_{\tilde{B}}
\preceq\left(\begin{array}{cc}
\mathbf{S}^{\rmH}\mathbf{S}&-(\mathbf{I}+\mathbf{S}^{\rmH}\mathbf{S})^{-1}\mathbf{S}^{\rmH}\\-\mathbf{S}(\mathbf{I}+\mathbf{S}^{\rmH}\mathbf{S})^{-1}&\mathbf{I}
\end{array}\right).
\end{eqnarray}
For the Frobenius norm, Lemma \ref{lem:pm}, \eqref{range} and \eqref{eqn6} we have
\begin{eqnarray*}\label{p1}
\|(\mathbf{I}-\mathbf{P}_{Q})\mathbf{A}\|_{\rmF}^{2} 
&=& \|\mathbf{U}^{\rmH}(\mathbf{I}-\mathbf{P}_{Q})\mathbf{U}\boldsymbol{\Sigma} \mathbf{V}^{\rmH}\|_{\rmF}^{2}\\
&=& \|(\mathbf{I}-\mathbf{P}_{U^{\rmH}Q})\boldsymbol{\Sigma} \mathbf{V}^{\rmH}\|_{\rmF}^{2}\\
&\leq&\|(\mathbf{I}-\mathbf{P}_{\tilde{B}})\boldsymbol{\Sigma}\|_{\rmF}^{2}\\
&=&\mathrm{tr}(\boldsymbol{\Sigma}^{\rmH}(\mathbf{I}-\mathbf{P}_{\tilde{B}})\boldsymbol{\Sigma})\\
&=&\mathrm{tr}(\boldsymbol{\Sigma}_{r}^{\rmH}\mathbf{S}^{\rmH}\mathbf{S}\boldsymbol{\Sigma}_{r})+\mathrm{tr}(\boldsymbol{\Sigma}_{0}^{\rmH}\boldsymbol{\Sigma}_{0})\\
&=&\|\mathbf{S}\boldsymbol{\Sigma}_{r}\|_{\rmF}^{2}+\|\boldsymbol{\Sigma}_{0}\|_{\rmF}^{2}\\
&\leq&\alpha^{2\tau}\|\boldsymbol{\Sigma}_{0}\mathbf{\Lambda}_{2}\mathbf{\Lambda}_{1}^{\dagger}\|_{\rmF}^{2}+\|\boldsymbol{\Sigma}_{0}\|_{\rmF}^{2}.
\end{eqnarray*}
where $\alpha=\frac{\sigma_{r+1}}{\sigma_{r}}$. Recall that $\mathbf{\Lambda}_{1} = \mathbf{V}_{r}^{\rmH}\boldsymbol{\Omega}$ and $\mathbf{\Lambda}_{2} = \mathbf{V}_{0}^{\rmH}\boldsymbol{\Omega}$. The Gaussian distribution is rotationally invariant, so $\mathbf{V}^{\rmH}\boldsymbol{\Omega}$ is also a standard Gaussian
matrix. Observe that $\mathbf{\Lambda}_{1}$ and $\mathbf{\Lambda}_{2}$ are nonoverlapping submatrices of $\mathbf{V}^{\rmH}\boldsymbol{\Omega}$, so
these two matrices are not only standard Gaussian but also stochastically independent.
By \eqref{th3.3:1} and Lemmas \ref{lem:sgt} and \ref{lem:ginv} we have
\begin{eqnarray*}
    \mathbb{E}\|\mathbf{A}-\tilde{\mathbf{A}}\|_{\rmF} &\leq& \sqrt{\alpha^{2\tau}\mathbb{E}\|\boldsymbol{\Sigma}_{0}\mathbf{\Lambda}_{2}\mathbf{\Lambda}_{1}^{\dagger}\|_{\rmF}^{2}+\|\boldsymbol{\Sigma}_{0}\|_{\rmF}^{2}}\\
    &=& \sqrt{\alpha^{2\tau}\|\boldsymbol{\Sigma}_{0}\|_{\rmF}^{2}\mathbb{E}\|\mathbf{\Lambda}_{1}^{\dagger}\|_{\rmF}^{2} + \|\boldsymbol{\Sigma}_{0}\|_{\rmF}^{2}}\\
    &=& \left(1+\frac{r\alpha^{2\tau}}{d-r-1}\right)^{\frac{1}{2}}\|\boldsymbol{\Sigma}_{0}\|_{\rmF}\\
    &=& \left(1+\frac{r\alpha^{2\tau}}{d-r-1}\right)^{\frac{1}{2}}\left(\sum_{j>r}\sigma_{j}^{2}\right)^{\frac{1}{2}}.
\end{eqnarray*}
This completes the proof.
\end{IEEEproof}

The following theorem will give the relationship between the error of singular values obtained by Algorithm \ref{alg:qr} and the error of low-rank approximation of matrix $\mathbf{A}$.

\begin{theorem}\label{singular}
Let singular values of $\mathbf{A}\in\mathbb{C}^{m\times n}\;(m\geq n)$ have decreasing orders with $\sigma_{1}(\mathbf{A})\geq\sigma_{2}(\mathbf{A})\geq\cdots\geq\sigma_{n}(\mathbf{A})$. Let $\boldsymbol{\Omega}\in\mathbb{R}^{n\times d}$ be standard Gaussian matrix. Let singular values of $\tilde{\mathbf{A}}=\mathbf{Q}\mathbf{Q}^H\mathbf{A}$ have decreasing orders with $\sigma_{1}(\tilde{\mathbf{A}})\geq\sigma_{2}(\tilde{\mathbf{A}})\geq\cdots\geq\sigma_{n}(\tilde{\mathbf{A}})$.
Then
$$
\frac{\max\limits_{1\leq i\leq d}| \sigma_{i}^2(\mathbf{A})-\sigma_{i}^2(\tilde{\mathbf{A}})|}{\| \mathbf{A}\|_{\rmF}^2} \leq\|\mathbf{Q}^{\bot}\|_2^2,
$$
where $\mathbf{Q}^{\bot}=\mathbf{I}-\mathbf{Q}\mathbf{Q}^H$ gradually approaches zero matrix as $d$ gradually increases to $n$.
\end{theorem}
\begin{IEEEproof}
Let $\mathbf{K}=(q_{1},q_{2},\ldots,q_{d},q_{d+1},\ldots,q_{n})$, $\mathbf{Q}=(q_{1},q_{2},\ldots,q_{d})$, $\mathbf{Q}_{0}=(q_{d+1},\ldots,q_{n})$.
Then
\begin{align*}
  \mathbf{K}\mathbf{K}^{\rmH}=(\mathbf{Q},\mathbf{Q}_{0})(\mathbf{Q},\mathbf{Q}_{0})^{\rmH}=\mathbf{Q}\mathbf{Q}^{\rmH}+\mathbf{Q}_{0}\mathbf{Q}_{0}^{\rmH}.
\end{align*}
Thus,
\begin{align*}
  \mathbf{K}\mathbf{K}^{\rmH}\succeq \mathbf{Q}\mathbf{Q}^{\rmH}.
\end{align*}
Since $\mathbf{A}=\mathbf{K}\mathbf{K}^{\rmH}\mathbf{A}$ and $\tilde{\mathbf{A}}=\mathbf{Q}\mathbf{Q}^{\rmH}\mathbf{A}$. Let $\mathbf{C}=\mathbf{K}^{\rmH}\mathbf{A}$, $\tilde{\mathbf{C}}=\mathbf{Q}^{\rmH}\mathbf{A}$, then for $i=1,2,\ldots,d$, we have
\begin{eqnarray*}
  \sigma_{i}^{2}(\mathbf{A})&=&\lambda_{i}(\mathbf{C}\mathbf{C}^{\rmH})\\
  &=&\lambda_{i}(\mathbf{K}^{\rmH}\mathbf{A}\mathbf{A}^{\rmH}\mathbf{K})\\
  &=&\lambda_{i}(\mathbf{A}^{\rmH}\mathbf{K}\mathbf{K}^{\rmH}\mathbf{A})
\end{eqnarray*}
and
\begin{eqnarray*}
  \sigma_{i}^{2}(\tilde{\mathbf{A}})&=&\lambda_{i}(\tilde{\mathbf{C}}\tilde{\mathbf{C}}^{\rmH})\\
  &=&\lambda_{i}(\mathbf{Q}^{\rmH}\mathbf{A}\mathbf{A}^{\rmH}\mathbf{Q})\\
  &=&\lambda_{i}(\mathbf{A}^{\rmH}\mathbf{Q}\mathbf{Q}^{\rmH}\mathbf{A}).
\end{eqnarray*}
Then
\begin{align*}
  \sigma_{i}^{2}(\mathbf{A})\geq \sigma_{i}^{2}(\tilde{\mathbf{A}}), i=1,2,\ldots,d.
\end{align*}
It follows that
\begin{eqnarray*}
   \|(\mathbf{I}-\mathbf{Q}\mathbf{Q}^{\rmH})\mathbf{A}\|_{\rmF}^2 & = & \left|\|\mathbf{A}\|_{\rmF}^2 - \|\tilde{\mathbf{A}}\|_{\rmF}^2\right| \\
   & = & \left|\sum\limits_{i=1}^{n}\sigma_{i}^{2}(\mathbf{A}) - \sum\limits_{i=1}^{d}\sigma_{i}^{2}(\tilde{\mathbf{A}})\right|\\
   & = & \sum\limits_{i=1}^{d}\left|\sigma_{i}^{2}(\mathbf{A})-\sigma_{i}^{2}(\tilde{\mathbf{A}})\right| + \sum\limits_{i=d+1}^{n}\sigma_{i}^{2}(\mathbf{A})\\
   &\geq & \sum\limits_{i=1}^{d}\left|\sigma_{i}^{2}(\mathbf{A})-\sigma_{i}^{2}(\tilde{\mathbf{A}})\right|\\
   &\geq & \max\limits_{1\leq i\leq d} \left|\sigma_{i}^{2}(\mathbf{A})-\sigma_{i}^{2}(\tilde{\mathbf{A}})\right|.
\end{eqnarray*}

Thus,
\begin{eqnarray*}
\max\limits_{1\leq i\leq d}| \sigma_{i}^2(\mathbf{A})-\sigma_{i}^2(\tilde{\mathbf{A}})|
&\leq&\|(\mathbf{I}-\mathbf{Q}\mathbf{Q}^{\rmH})\mathbf{A}\|_{\rmF}^2\\
&\leq&\|\mathbf{I}-\mathbf{Q}\mathbf{Q}^{\rmH}\|_{2}^{2} \|\mathbf{A}\|_{\rmF}^{2},
\end{eqnarray*}
which yields
\[\frac{\max\limits_{1\leq i\leq d}| \sigma_{i}^2(\mathbf{A})-\sigma_{i}^2(\tilde{\mathbf{A}})|}{\| \mathbf{A}\|_{\rmF}^2} \leq \|\mathbf{Q}^{\bot}\|_2^2,
\]
where $\mathbf{Q}^{\bot}=\mathbf{I}-\mathbf{Q}\mathbf{Q}^{\rmH}$ gradually approaches zero matrix as $d$ gradually increases to $n$.
This completes the proof.
\end{IEEEproof}

\section{Numerical experiments}


In this section, we present the results of some numerical experiments conducted to evaluate the performance of Algorithm \ref{alg:A} for approximating a low-rank input matrix. 
We show that Algorithm \ref{alg:A} provides highly accurate singular values and low-rank approximations, and compare our algorithm against several other algorithms from the literature. We furthermore employ Algorithm \ref{alg:A} for solving the image reconstruction. The experiments were run in MATLAB R2021b on a desktop PC with a 3.30 GHz AMD Ryzen 9 5900HX processor and 16 GB of memory. The source code of our method is published at \url{https://github.com/xuweiwei1/EOD-ABE.git}.

\subsection{Test Matrices With Randomly Generated Variables}
For the sake of simplicity, we focus on square matrices. We construct four classes of input matrices in \cite{kc21} to illustrate the suitability and robustness of Algorithm \ref{alg:A}. The first two classes contain one or multiple gaps in the spectrum and are particularly designed to investigate the rank-revealing property of Algorithm \ref{alg:A}. The second two classes have fast and slow decay singular values. We generate square matrices of order $n = 1000$.
\begin{itemize}
  \item Matrix 1 (Low-rank plus noise). This rank-$r$ matrix, with $r=20$, is formed as follows:
  \begin{equation}\label{matrix-1}
    \mathbf{A} = \mathbf{U}\boldsymbol{\Sigma} \mathbf{V}^{\rmT}+\alpha\sigma_{r} \mathbf{E},
  \end{equation}
  where $\mathbf{U}\in \mathbb{R}^{n\times n}$ and $\mathbf{V} \in \mathbb{R}^{n\times n}$ are orthogonal matrices, and $\boldsymbol{\Sigma}=\mathrm{diag}\{\sigma_{1},\ldots,\sigma_{n}\} \in \mathbb{R}^{n\times n}$ is diagonal whose entries ($\sigma_{i}$s) decrease
linearly from 1 to $10^{-25}$, $\sigma_{r+1}=...=\sigma_{n}=0$, and $\mathbf{E}\in\mathbb{R}^{n\times n}$ is standard Gaussian matrix.
  \begin{enumerate}
    \item $\alpha=0.005$ in which the matrix has a gap $\approx 200$.
    \item $\alpha=0.02$ in which the matrix has a gap $\approx 50$.
  \end{enumerate}
  \item Matrix 2 (The devil's stairs). This challenging matrix has multiple gaps in its spectrum. The singular values are arranged analogues to a descending staircase with each step consisting of 15 equal singular values.
  \item Matrix 3 (Fast decay). This matrix is formed as follows:
  \begin{equation}\label{matrix-3}
    \mathbf{A} = \mathbf{U}\boldsymbol{\Sigma} \mathbf{V}^{\rmT},
  \end{equation}
  where the diagonal elements of $\boldsymbol{\Sigma}$ have the form $\sigma_{i}=e^{-i/6}$, for $i=1,...,n$.
  \item Matrix 4 (Slow decay). This matrix is also formed as Matrix 3, but the diagonal elements of $\boldsymbol{\Sigma}$ take the form $\sigma_{i}=i^{-2}$, for $i=1,...,n$.
\end{itemize}

%

\begin{figure}[!ht]
  \centering
  \includegraphics[width=1.7in]{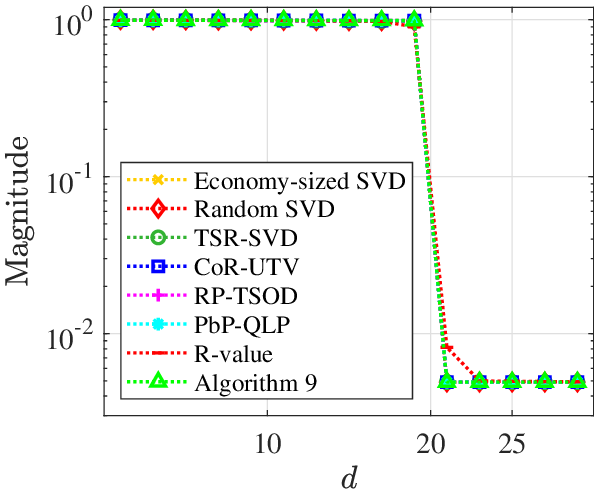}
  \includegraphics[width=1.7in]{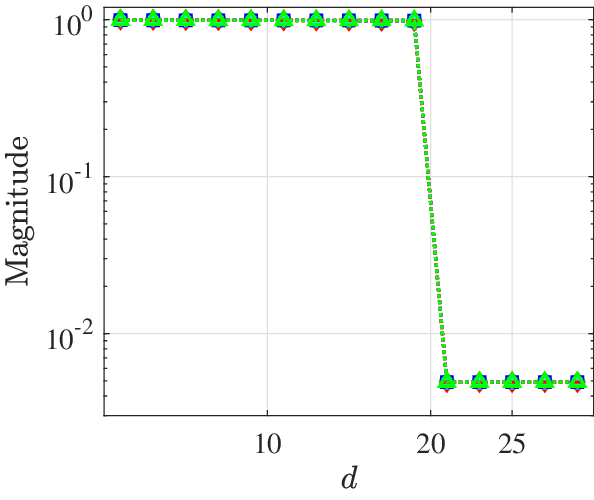}\\
  \caption{Singular value approximations for Matrix 1 with $\alpha=0.005$. Left: $\tau=0$. Right: $\tau=2$.}\label{plus_noiseLargeGap}
\end{figure}
\begin{figure}[!ht]
  \centering
  \includegraphics[width=1.7in]{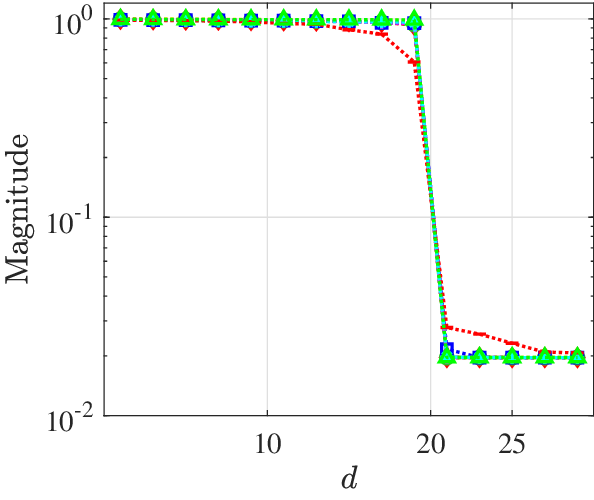}
  \includegraphics[width=1.7in]{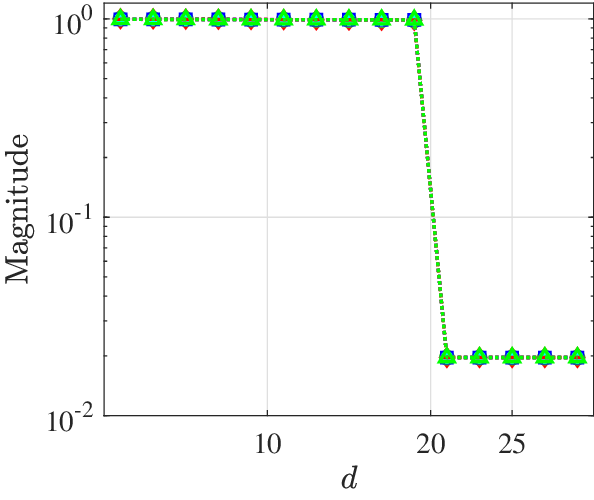}\\
  \caption{Singular value approximations for Matrix 1 with $\alpha=0.02$. Left: $\tau=0$. Right: $\tau=2$.}\label{plus_noiseSmallGap}
\end{figure}
\begin{figure}[!ht]
  \centering
  \includegraphics[width=1.7in]{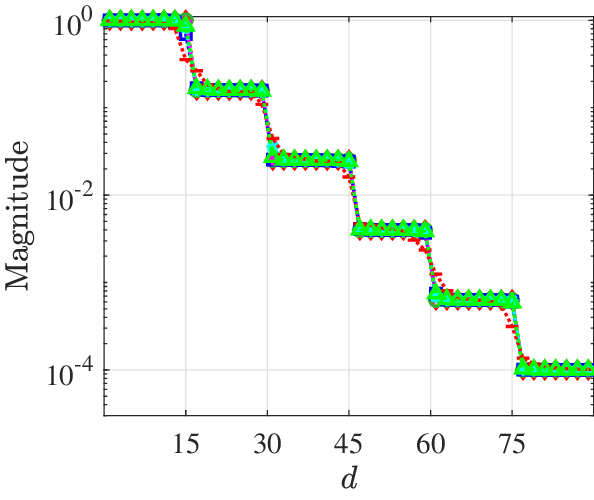}
  \includegraphics[width=1.7in]{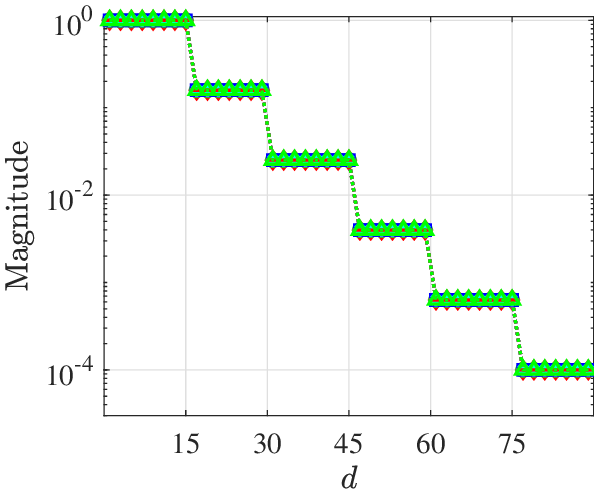}\\
  \caption{Singular value approximations for Matrix 2. Left: $\tau=0$. Right: $\tau=2$.}\label{stairs}
\end{figure}
\begin{figure}[!ht]
  \centering
  \includegraphics[width=1.7in]{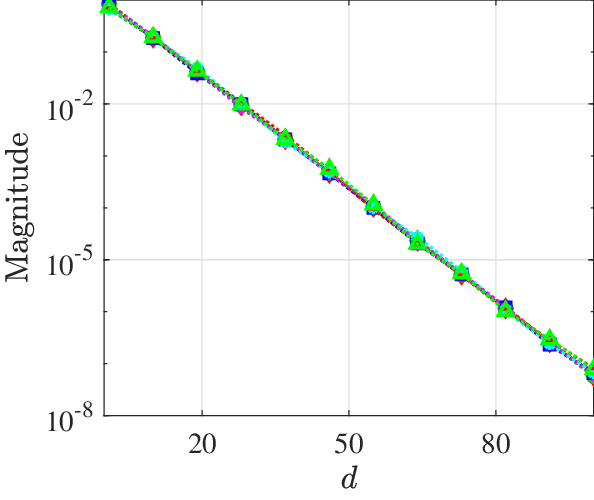}
  \includegraphics[width=1.7in]{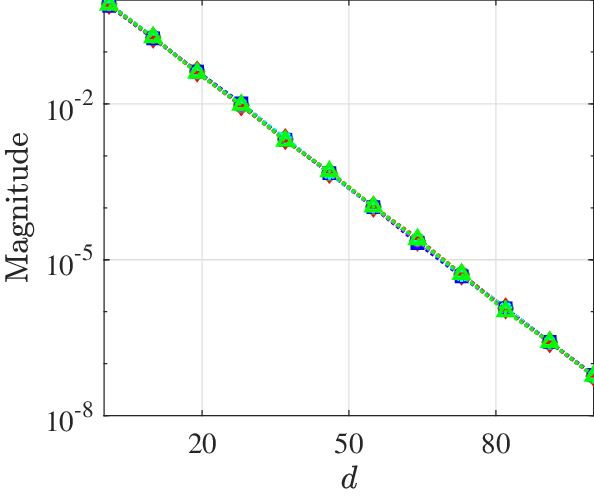}\\
  \caption{Singular value approximations for Matrix 3. Left: $\tau=0$. Right: $\tau=2$.}\label{fast_decay}
\end{figure}
\begin{figure}[!ht]
  \centering
  \includegraphics[width=1.7in]{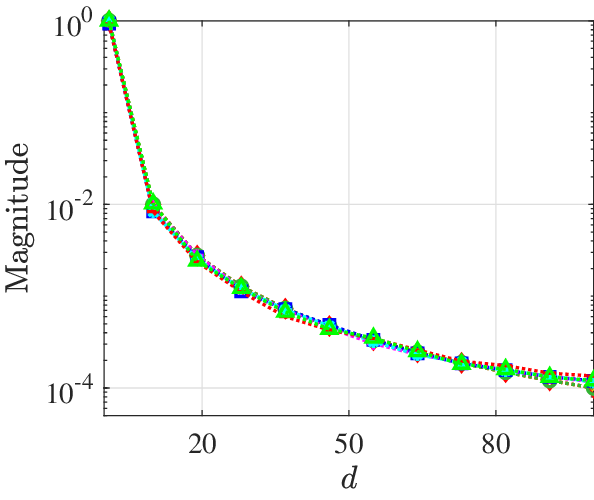}
  \includegraphics[width=1.7in]{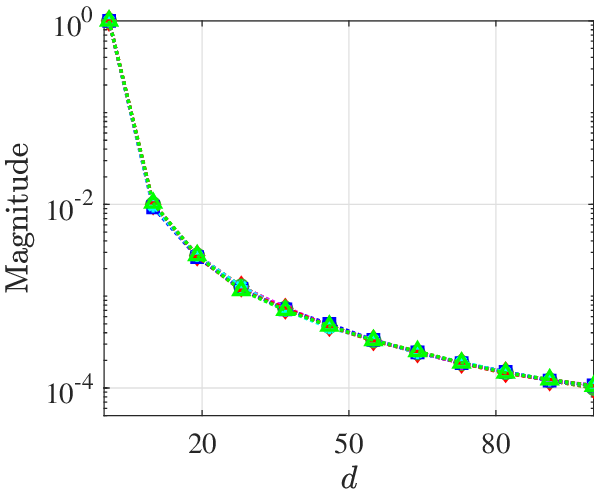}\\
  \caption{Singular value approximations for Matrix 4. Left: $\tau=0$. Right: $\tau=2$.}\label{slow_decay}
\end{figure}

The results for singular values estimation are plotted in Figs. \ref{plus_noiseLargeGap}, \ref{plus_noiseSmallGap}, \ref{stairs}, \ref{fast_decay} and \ref{slow_decay}. We make several observations:

\begin{enumerate}
  \item The numerical rank of both matrices 1 for $\alpha=0.005$ and $\alpha=0.02$ is strongly revealed in $\boldsymbol{\Sigma}$ generated by Algorithm \ref{alg:A} with $\tau=0$. This is due to the fact that the gaps in the spectrums of these matrices are well-defined. The diagonal elements of $\boldsymbol{\Sigma}$ generated by Algorithm \ref{alg:A} with subspace iteration technique are as accurate as those of the optimal SVD. Figs. \ref{plus_noiseLargeGap} and \ref{plus_noiseSmallGap} show that Algorithm \ref{alg:A} is a rank-revealer.
  \item For Matrix 2, R-values of Algorithm \ref{alg:A} with $\tau=0$ do not clearly disclose the gaps in matrix’s spectrum. This is because the gaps are not substantial. However, Algorithm \ref{alg:A} with $\tau=0$ strongly reveals the gaps, which shows that the procedure leading to the formation of the upper triangular matrix $\mathbf{R}$ provides a good first step for Algorithm \ref{alg:A}. R-values of Algorithm \ref{alg:A} with subspace iteration clearly disclose that gaps.
  \item For Matrices 3 and 4, Algorithm \ref{alg:A} provides highly accurate singular values, showing similar performance as random SVD in \cite{hmt11}, TSR-SVD in \cite{hmt11}, CoR-UTV in \cite{kl18}, RP-TSOD in \cite{kc20}, PbP-QLP in \cite{kc21} and the built-in function economy-sized SVD (\texttt{[U,S,V]=svd(A,'econ')}) in MATLAB.
\end{enumerate}

The results presented in Figs. \ref{plus_noiseLargeGap}, \ref{plus_noiseSmallGap}, \ref{stairs}, \ref{fast_decay} and \ref{slow_decay} demonstrate the applicability of Algorithm \ref{alg:A} in accurately estimating the singular values of matrices from different classes.

\subsection{Comparison of low-rank approximation with existing methods}
We compare the speed and accuracy of Algorithm \ref{alg:A} against the existing state-of-the-art algorithms, in factoring input matrices with various dimensions.
We construct a new class of strictly rank-deficient matrices
$\mathbf{A} \in \mathbb{R}^{n\times n}$ with rank $r$ of the form:
\begin{eqnarray}\label{mat2}
\mathbf{A} = \mathbf{U}\left(
         \begin{array}{cc}
           \boldsymbol{\Sigma} & \mathbf{O}_{r\times (n-r)} \\
           \mathbf{O}_{(n-r)\times r} & \mathbf{O}_{(n-r)\times (n-r)} \\
         \end{array}
       \right)
 \mathbf{V}^{\rmT},
\end{eqnarray}
where $\mathbf{U}\in \mathbb{R}^{n\times n}$ and $\mathbf{V} \in \mathbb{R}^{n\times n}$ are orthogonal matrices, and $\boldsymbol{\Sigma}=\mathrm{diag}\{\sigma_{1},\ldots,\sigma_{r}\} \in \mathbb{R}^{r\times r}$ is diagonal matrix, $\sigma_{1}\geq\sigma_{2}\geq\cdots\geq\sigma_{r}>0$ are the $r$ non-zero singular values of $\mathbf{A}$.
We generate the exact $\boldsymbol{\Sigma}=\mathrm{diag}\{\sigma_{1},\ldots,\sigma_{r}\}$ by using the built-in function \texttt{rand(1,r)}
in MATLAB such that
\[\sigma_{1}\geq\sigma_{2}\geq\cdots\geq\sigma_{r}>0\]
and the orthogonal matrices $\mathbf{U}$ and $\mathbf{V}$ are generated by code \texttt{orth(randn(n,n))} in MATLAB.


We aim to assess the following issues for matrix \eqref{mat2} of different orders.

\begin{itemize}
  \item Assessing calculation time (second) for computing low-rank approximation by Algorithm \ref{alg:A} and the algorithms mentioned earlier.
  \item Define the relative error value as follows:
  \[\mathrm{RelErr}:=\frac{\|\mathbf{A}-\tilde{\mathbf{A}}\|_{\rmF}}{\|\mathbf{A}\|_{\rmF}},\]
  where $\tilde{\mathbf{A}}$ is rank-$r$ approximations of $\mathbf{A}$ by either method.
  Assessing the relative error values of the above eight algorithms.
\end{itemize}



First, we consider that when the rank of the matrix is unknown, the sampling size parameter $d$ may have the following three cases: 
(1) when the sampling size parameter equals the matrix rank, i.e., $d = r$;
(2) when the sampling size parameter exceeds the matrix rank, i.e., $d > r$; and
(3) when the sampling size parameter is less than the matrix rank, i.e., $d < r$.
To evaluate these scenarios, we used square matrices of various sizes, setting the matrix rank at $r = 0.4n$.
Sampling size parameters of $d = 0.4n$, $d = 0.6n$, and $d = 0.35n$ are tested.
The computation time and accuracy of Algorithm \ref{alg:A} are compared with the current state-of-the-art algorithms from the literature.
The results are presented in Tables \ref{0.4n}, \ref{0.6n} and \ref{0.35n}.


\begin{table}
\setlength\tabcolsep{3pt}
  \centering
  \caption{Numerical results of low-rank approximation with different methods for matrix \eqref{mat2} ($r=0.4n$, $d=0.4n$)}\label{0.4n}
  \begin{tabular}{lccccccc}
    \hline
    \multirow{2}{1.4cm}{Algorithms}  & & \multicolumn{2}{c}{$n=4000$} & \multicolumn{2}{c}{$n=8000$} & \multicolumn{2}{c}{$n=12000$} \\
    \cline{3-8}
     & & Time & RelErr& Time & RelErr& Time & RelErr  \\
    \hline
    \multirow{2}{1.4cm}{Economy-sized SVD} &\multirow{2}{0.7cm}{} & \multirow{2}{*}{9.9} & \multirow{2}{*}{5.8E-15} & \multirow{2}{*}{70.5} & \multirow{2}{*}{7.4E-15}& \multirow{2}{*}{235.9} & \multirow{2}{*}{8.5E-15}\\
                                    &          &       &       &   &      &&\\
    \multirow{3}{1.4cm}{Random SVD}	&	$\tau=0$	&	2.2 	&	1.9E-14	&	26.0 	&	3.6E-14	&	98.4 	&	5.0E-14	\\
	&	$\tau=1$	&	4.1 	&	4.7E-15	&	50.4 	&	5.8E-15	&	194.6 	&	6.8E-15	\\
	&	$\tau=2$	&	6.1 	&	4.7E-15	&	75.6 	&	5.7E-15	&	289.1 	&	6.7E-15	\\
	&		&		&		&		&		&		&		\\
\multirow{3}{1.4cm}{TSR-SVD}	&	$\tau=0$	&	3.1 	&	3.8E-14	&	38.6 	&	4.8E-14	&	126.3 	&	2.7E-14	\\
	&	$\tau=1$	&	5.1 	&	7.4E-15	&	70.0 	&	9.3E-15	&	232.4 	&	1.3E-14	\\
	&	$\tau=2$	&	7.0 	&	7.3E-15	&	90.0 	&	9.3E-15	&	341.3 	&	1.3E-14	\\
	&		&		&		&		&		&		&		\\
\multirow{3}{1.4cm}{CoR-UTV}	&	$\tau=0$	&	2.8 	&	1.9E-12	&	31.4 	&	6.6E-12	&	124.8 	&	1.6E-12	\\
	&	$\tau=1$	&	4.7 	&	3.1E-14	&	56.9 	&	4.3E-14	&	216.1 	&	5.4E-14	\\
	&	$\tau=2$	&	6.8 	&	5.9E-15	&	85.4 	&	7.3E-15	&	306.9 	&	8.6E-15	\\
	&		&		&		&		&		&		&		\\
\multirow{3}{1.4cm}{RP-TSOD}	&	$\tau=0$	&	2.2 	&	1.8E-14	&	20.2 	&	3.8E-14	&	74.2 	&	7.7E-14	\\
	&	$\tau=1$	&	3.1 	&	1.8E-15	&	27.3 	&	2.1E-15	&	97.7 	&	2.4E-15	\\
	&	$\tau=2$	&	3.9 	&	1.7E-15	&	34.6 	&	2.0E-15	&	121.3 	&	2.3E-15	\\
	&		&		&		&		&		&		&		\\
\multirow{3}{1.4cm}{PbP-QLP}	&	$\tau=0$	&	1.7 	&	1.8E-14	&	18.6 	&	4.0E-14	&	62.2 	&	3.9E-14	\\
	&	$\tau=1$	&	3.6 	&	4.2E-15	&	43.0 	&	5.2E-15	&	158.0 	&	6.1E-15	\\
	&	$\tau=2$	&	5.7 	&	4.3E-15	&	74.9 	&	5.2E-15	&	243.8 	&	6.0E-15	\\
	&		&		&		&		&		&		&		\\
\multirow{3}{1.4cm}{Algorithm \ref{alg:A}}	&	$\tau=0$	&	\textbf{1.2} 	&	3.1E-13	&	\textbf{8.9} 	&	1.8E-12	&	\textbf{25.4} 	&	1.7E-12	\\
	&	$\tau=1$	&	2.0 	&	1.3E-15	&	16.4 	&	1.3E-15	&	46.3 	&	1.3E-15	\\
	&	$\tau=2$	&	3.0 	&	1.2E-15	&	23.2 	&	1.3E-15	&	71.9 	&	1.3E-15	\\
    \hline
  \end{tabular}
\end{table}

\begin{table}
\setlength\tabcolsep{3pt}
  \centering
  \caption{Numerical results of low-rank approximation with different methods for matrix \eqref{mat2} ($r=0.4n$, $d=0.6n$)}\label{0.6n}
  \begin{tabular}{lcrrrrrr}
    \hline
    \multirow{2}{1.4cm}{Algorithms}  & & \multicolumn{2}{c}{$n=4000$} & \multicolumn{2}{c}{$n=8000$} & \multicolumn{2}{c}{$n=12000$} \\
    \cline{3-8}
     & & Time & RelErr& Time & RelErr& Time & RelErr  \\
    \hline
    \multirow{2}{1.4cm}{Economy-sized SVD} &\multirow{2}{0.7cm}{} & \multirow{2}{*}{11.3} & \multirow{2}{*}{5.8E-15} & \multirow{2}{*}{82.5} & \multirow{2}{*}{7.6E-15}& \multirow{2}{*}{232.0} & \multirow{2}{*}{8.6E-15}\\
                                    &          &       &       &   &      &&\\
    \multirow{3}{1.4cm}{Random SVD}	&	$\tau=0$	&	5.5 	&	5.9E-15	&	60.3 	&	7.4E-15	&	179.4 	&	8.7E-15	\\
	&	$\tau=1$	&	7.5 	&	4.7E-15	&	84.6 	&	5.8E-15	&	270.9 	&	6.7E-15	\\
	&	$\tau=2$	&	9.4 	&	4.7E-15	&	118.8 	&	5.8E-15	&	369.8 	&	6.9E-15	\\
	&		&		&		&		&		&		&		\\
\multirow{3}{1.4cm}{TSR-SVD}	&	$\tau=0$	&	6.0 	&	7.1E-15	&	67.4 	&	9.1E-15	&	221.1 	&	1.3E-14	\\
	&	$\tau=1$	&	8.4 	&	7.4E-15	&	90.7 	&	9.3E-15	&	322.0 	&	1.3E-14	\\
	&	$\tau=2$	&	9.8 	&	7.3E-15	&	113.6 	&	9.4E-15	&	401.0 	&	1.3E-14	\\
	&		&		&		&		&		&		&		\\
\multirow{3}{1.4cm}{CoR-UTV}	&	$\tau=0$	&	8.1 	&	2.3E-14	&	84.7 	&	3.5E-14	&	297.7 	&	4.5E-14	\\
	&	$\tau=1$	&	12.6 	&	5.8E-15	&	140.1 	&	7.2E-15	&	465.9 	&	8.4E-15	\\
	&	$\tau=2$	&	15.9 	&	5.7E-15	&	173.0 	&	7.4E-15	&	562.8 	&	8.4E-15	\\
	&		&		&		&		&		&		&		\\
\multirow{3}{1.4cm}{RP-TSOD}	&	$\tau=0$	&	5.3 	&	2.0E-15	&	45.5 	&	2.3E-15	&	152.0 	&	2.5E-15	\\
	&	$\tau=1$	&	6.9 	&	1.8E-15	&	55.9 	&	2.1E-15	&	186.0 	&	2.4E-15	\\
	&	$\tau=2$	&	9.5 	&	1.7E-15	&	68.2 	&	2.0E-15	&	228.7 	&	2.3E-15	\\
	&		&		&		&		&		&		&		\\
\multirow{3}{1.4cm}{PbP-QLP}	&	$\tau=0$	&	4.8 	&	3.9E-15	&	44.8 	&	4.8E-15	&	149.7 	&	5.6E-15	\\
	&	$\tau=1$	&	6.9 	&	4.3E-15	&	69.5 	&	5.3E-15	&	247.4 	&	6.0E-15	\\
	&	$\tau=2$	&	8.5 	&	4.1E-15	&	98.3 	&	5.4E-15	&	332.9 	&	6.1E-15	\\
	&		&		&		&		&		&		&		\\
\multirow{3}{1.4cm}{Algorithm \ref{alg:A}}	&	$\tau=0$	&	\textbf{1.3} 	&	3.9E-13	&	\textbf{8.7} 	&	2.1E-12	&	\textbf{27.6} 	&	3.7E-12	\\
	&	$\tau=1$	&	2.5 	&	1.3E-15	&	15.4 	&	1.3E-15	&	51.1 	&	1.3E-15	\\
	&	$\tau=2$	&	3.1 	&	1.2E-15	&	22.2 	&	1.3E-15	&	67.8 	&	1.3E-15	\\
    \hline
  \end{tabular}
\end{table}

\begin{table}
\setlength\tabcolsep{3pt}
  \centering
  \caption{Numerical results of low-rank approximation with different methods for matrix \eqref{mat2} ($r=0.4n$, $d=0.35n$)}\label{0.35n}
  \begin{tabular}{lcrrrrrr}
    \hline
    \multirow{2}{1.4cm}{Algorithms}  & & \multicolumn{2}{c}{$n=4000$} & \multicolumn{2}{c}{$n=8000$} & \multicolumn{2}{c}{$n=12000$} \\
    \cline{3-8}
     & & Time & RelErr& Time & RelErr& Time & RelErr  \\
    \hline
    \multirow{2}{1.4cm}{Economy-sized SVD} &\multirow{2}{0.7cm}{} & \multirow{2}{*}{10.5} & \multirow{2}{*}{5.8E-15} & \multirow{2}{*}{72.3} & \multirow{2}{*}{7.4E-15}& \multirow{2}{*}{247.3} & \multirow{2}{*}{8.4E-15}\\
                                    &          &       &       &   &      &&\\
    \multirow{3}{1.4cm}{Random SVD}	&	$\tau=0$	&	1.8 	&	1.4E-01	&	18.3 	&	1.4E-01	&	71.2 	&	1.4E-01	\\
	&	$\tau=1$	&	3.3 	&	5.1E-02	&	35.8 	&	5.2E-02	&	138.3 	&	5.2E-02	\\
	&	$\tau=2$	&	4.9 	&	4.6E-02	&	54.2 	&	4.6E-02	&	210.5 	&	4.6E-02	\\
	&		&		&		&		&		&		&		\\
\multirow{3}{1.4cm}{TSR-SVD}	&	$\tau=0$	&	2.5 	&	1.4E-01	&	26.4 	&	1.4E-01	&	92.7 	&	1.4E-01	\\
	&	$\tau=1$	&	4.1 	&	5.2E-02	&	44.1 	&	5.2E-02	&	160.9 	&	5.2E-02	\\
	&	$\tau=2$	&	5.8 	&	4.6E-02	&	62.0 	&	4.6E-02	&	239.9 	&	4.7E-02	\\
	&		&		&		&		&		&		&		\\
\multirow{3}{1.4cm}{CoR-UTV}	&	$\tau=0$	&	2.2 	&	1.4E-01	&	22.3 	&	1.4E-01	&	83.5 	&	1.4E-01	\\
	&	$\tau=1$	&	3.9 	&	1.4E-01	&	40.3 	&	1.4E-01	&	150.4 	&	1.4E-01	\\
	&	$\tau=2$	&	5.5 	&	5.2E-02	&	57.9 	&	5.2E-02	&	223.6 	&	5.2E-02	\\
	&		&		&		&		&		&		&		\\
\multirow{3}{1.4cm}{RP-TSOD}	&	$\tau=0$	&	1.8 	&	1.4E-01	&	15.9 	&	1.4E-01	&	56.3 	&	1.4E-01	\\
	&	$\tau=1$	&	2.8 	&	5.2E-02	&	22.2 	&	5.1E-02	&	72.3 	&	5.2E-02	\\
	&	$\tau=2$	&	3.6 	&	4.7E-02	&	28.4 	&	4.6E-02	&	88.2 	&	4.7E-02	\\
	&		&		&		&		&		&		&		\\
\multirow{3}{1.4cm}{PbP-QLP}	&	$\tau=0$	&	1.4 	&	1.4E-01	&	13.1 	&	1.4E-01	&	43.7 	&	1.4E-01	\\
	&	$\tau=1$	&	3.1 	&	5.2E-02	&	30.9 	&	5.2E-02	&	107.6 	&	5.2E-02	\\
	&	$\tau=2$	&	4.6 	&	4.6E-02	&	49.5 	&	4.6E-02	&	170.4 	&	4.7E-02	\\
	&		&		&		&		&		&		&		\\
\multirow{3}{1.4cm}{Algorithm \ref{alg:A}}	&	$\tau=0$	&	\textbf{1.4} 	&	7.2E-13	&	\textbf{9.3} 	&	1.1E-12	&	\textbf{26.0} 	&	3.1E-12	\\
	&	$\tau=1$	&	2.3 	&	1.3E-15	&	17.6 	&	1.3E-15	&	47.2 	&	1.3E-15	\\
	&	$\tau=2$	&	3.3 	&	1.2E-15	&	26.2 	&	1.3E-15	&	68.9 	&	1.3E-15	\\
    \hline
  \end{tabular}
\end{table}

Seen from the Tables \ref{0.4n}, \ref{0.6n} and \ref{0.35n} we conclude the following results.

\begin{itemize}
  \item The algorithms mentioned earlier can all be calculated when the rank is known, while Algorithm \ref{alg:A} can calculate the singular value decomposition when the rank is unknown. 
  \item It can be observed that when the sampling size parameter is greater than or equal to the rank of the matrix, all algorithms achieve high computation accuracy. However, when the sampling size parameter is less than the rank of the matrix, the computation accuracy of random SVD, TSR-SVD, CoR-UTV, RP-TSOD, and PbP-QLP is significantly reduced, whereas Algorithm \ref{alg:A} maintains high accuracy. This is because Algorithm 9 can independently determine the position of the rank and stop computation at that point, avoiding undersampling.
  \item Algorithm \ref{alg:A} has a more obvious advantage in time. While ensuring high accuracy, it can be observed that in some cases it is 5x to 6x faster than the second-fastest algorithm. In the least case, algorithm \ref{alg:A} saves about 30\% runtime. As shown in Table \ref{0.35n}, Algorithm \ref{alg:A} achieves high accuracy while also providing faster computation times.
\end{itemize}

\begin{figure}[!ht]
  \centering
  \includegraphics[width=1.7in]{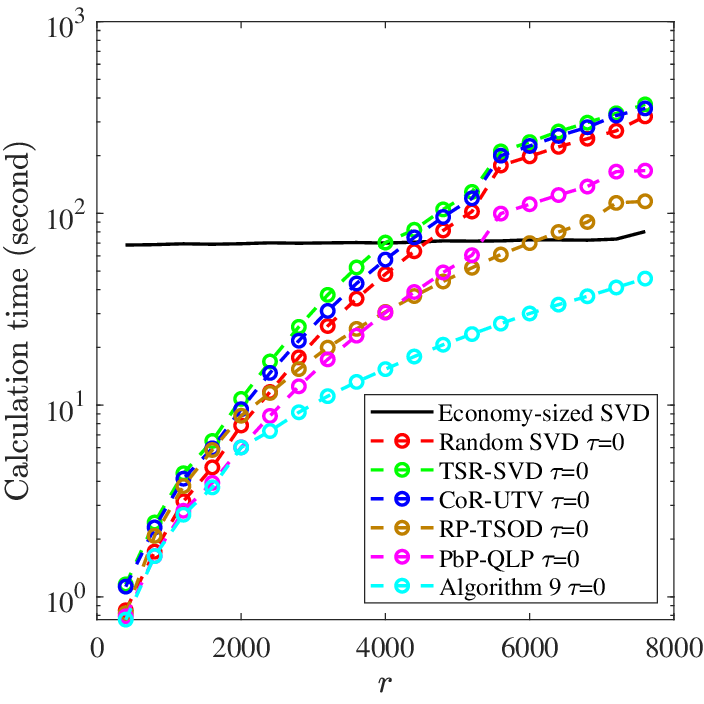} \includegraphics[width=1.7in]{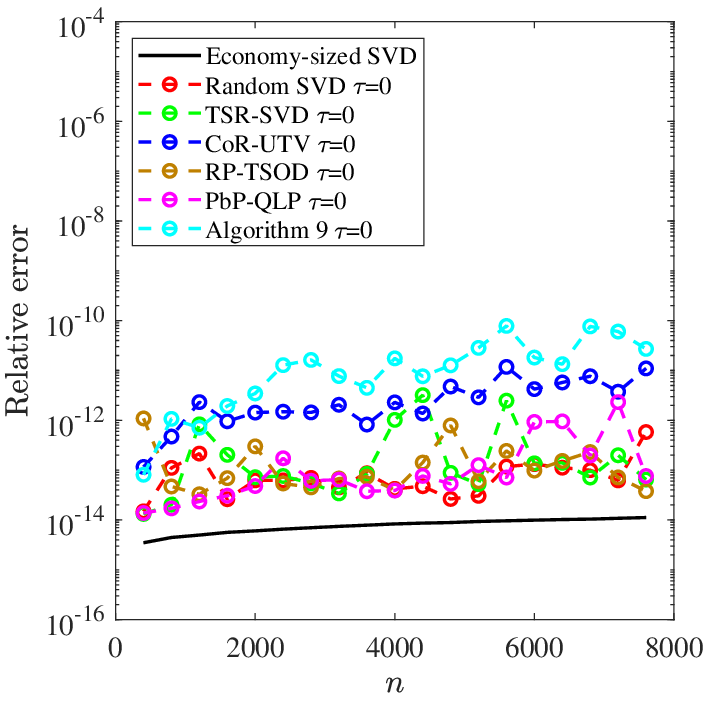}\\
  \caption{For $n=8000$, $\tau=0$, as the rank $r$ increases, the calculation time of different algorithms for matrix \eqref{mat2}}\label{8ct-1}
\end{figure}

\begin{figure}[!ht]
  \centering
  \includegraphics[width=1.7in]{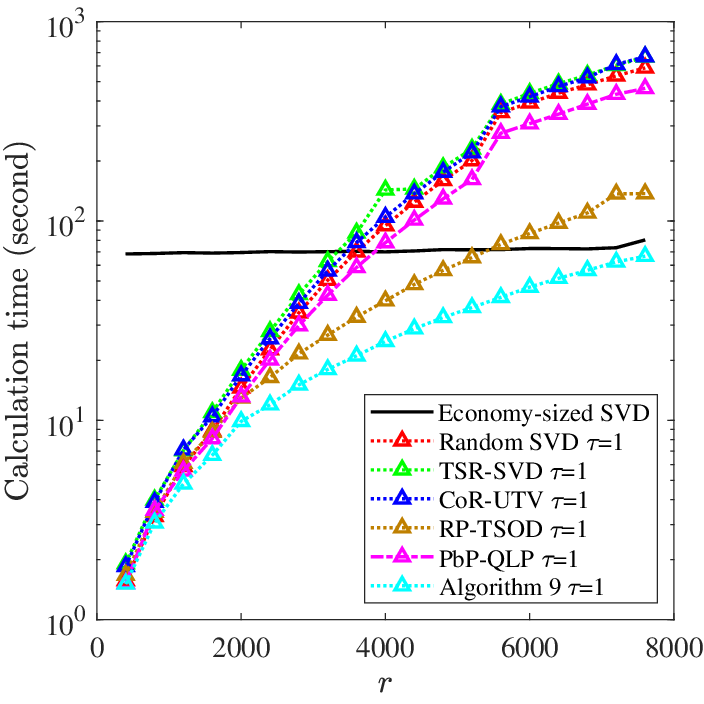} \includegraphics[width=1.7in]{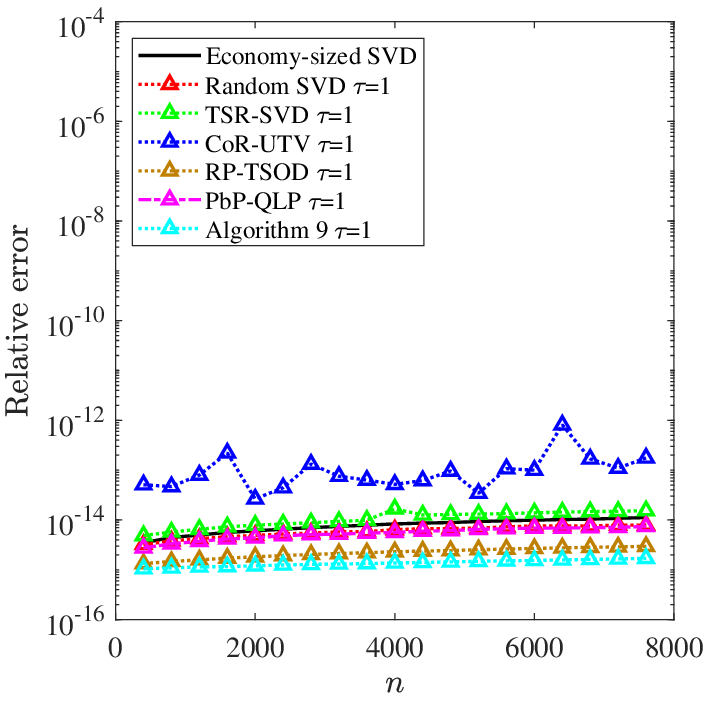}\\
  \caption{For $n=8000$, $\tau=1$, as the rank $r$ increases, the calculation time of different algorithms for matrix \eqref{mat2}}\label{8ct-2}
\end{figure}

\begin{figure}[!ht]
  \centering
  \includegraphics[width=1.7in]{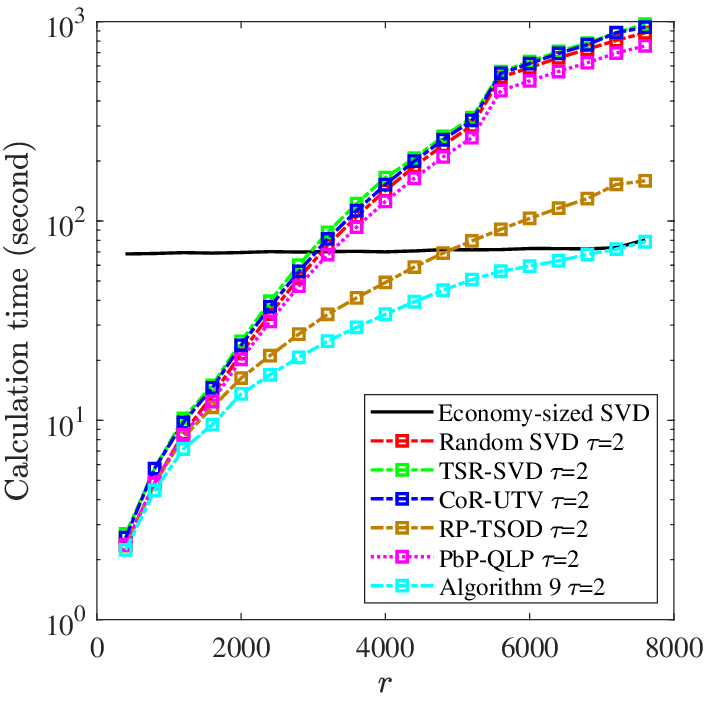} \includegraphics[width=1.7in]{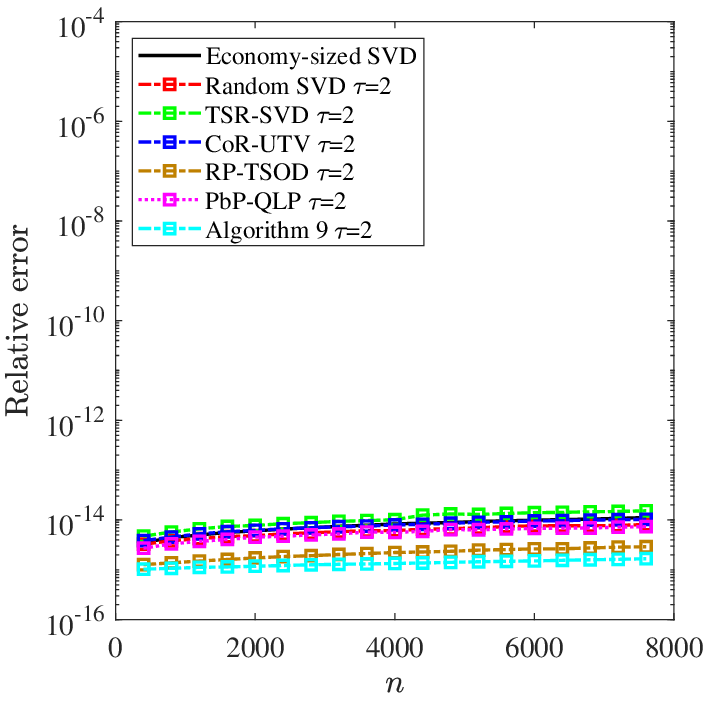}\\
  \caption{For $n=8000$, $\tau=2$, as the rank $r$ increases, the calculation time of different algorithms for matrix \eqref{mat2}}\label{8ct-3}
\end{figure}

Next, We investigate the changes in algorithm computation time and accuracy as the rank of the matrix varied.
We set $n = 8000$ and increase the rank from 400 to $n$ in increments of 400.
At the same time, the sampling parameter $d$ is set equal to the rank $r$ to calculate the matrix approximation.
Experiments are conducted using subspace iteration counts of $\tau = 0$, $\tau = 1$, and $\tau = 2$. The comparison results of computation time and accuracy are presented in Figures \ref{8ct-1}, \ref{8ct-2} and \ref{8ct-3}.

Seen from the Figures \ref{8ct-1}, \ref{8ct-2} and \ref{8ct-3} we conclude the following results.

\begin{itemize}
  \item It can be observed that Algorithm \ref{alg:A} takes less computing time than several other algorithms from the literature. As the matrix rank $r$ increases, the runtime gap between Algorithm \ref{alg:A} and other algorithms further enlarges. Therefore, the advantage of Algorithm \ref{alg:A} in computing time is more obvious when the rank is unknown.
  \item When computing strict low-rank approximations of matrices, relative error values are all small by all mentioned algorithms. In the worst case, Algorithm \ref{alg:A} achieves $10^{-12}$ absolute accuracy comparing to $10^{-15}$ of the best. Such an accuracy is sufficient in almost all applications. This means the calculation precisions by Algorithm \ref{alg:A} and several other algorithms from the literature are all high.
  \item Algorithm \ref{alg:A} with $\tau = 2$ gives highly accurate approximations.
\end{itemize}

Therefore, through the above numerical analysis, Algorithm \ref{alg:A} has better performance than other numerical methods when the rank is unknown.

\subsection{Application: Image Reconstruction}
\begin{figure}
  \centering
  \includegraphics[width=2.5in]{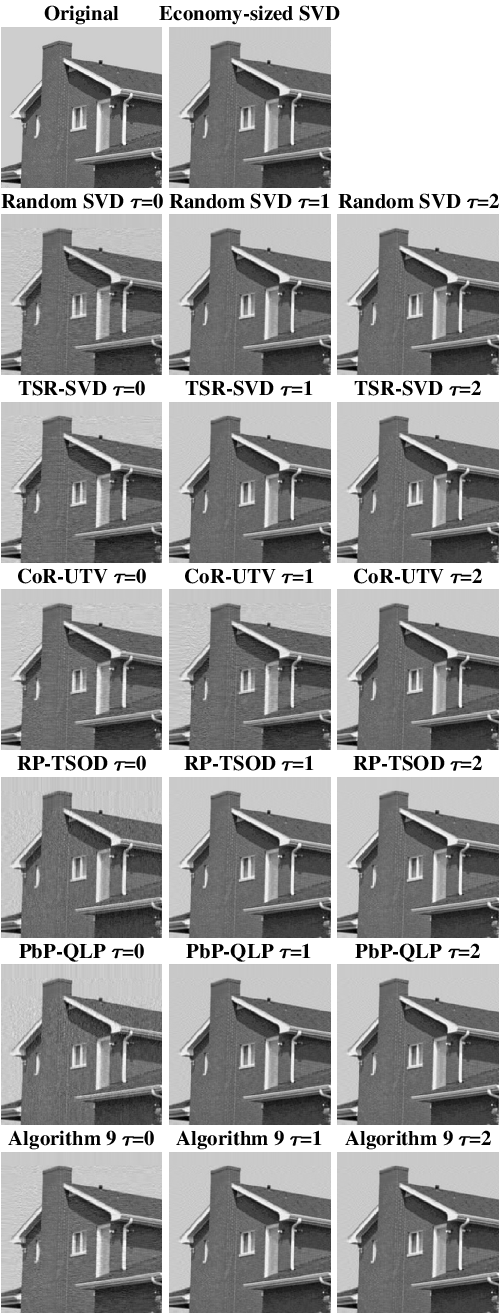}\;\;
  \caption{Low-rank image reconstruction. These figures show the results of
reconstructing a house image with dimension $512\times 512$ using $\varepsilon=0.001$ by different methods.}\label{example-house}
\end{figure}

\begin{figure}
  \centering
  \includegraphics[width=2.5in]{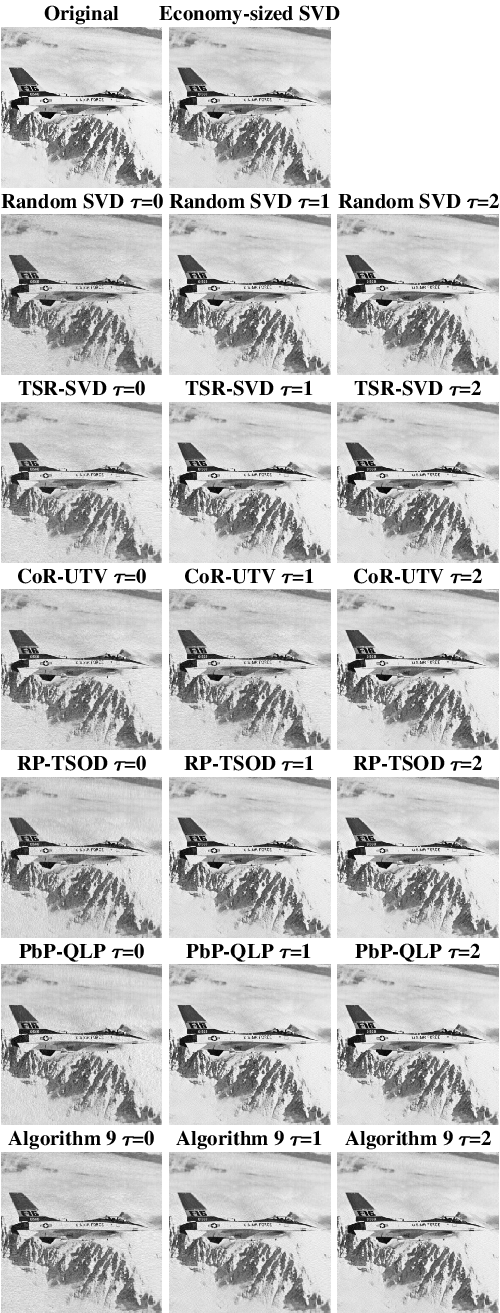}\;\;
  \caption{Low-rank image reconstruction. These figures show the results of
reconstructing a jetplane image with dimension $512\times 512$ using $\varepsilon=0.001$ by different methods.}\label{example-jetplane}
\end{figure}

\begin{figure}[!ht]
  \centering
  \includegraphics[width=3.3in]{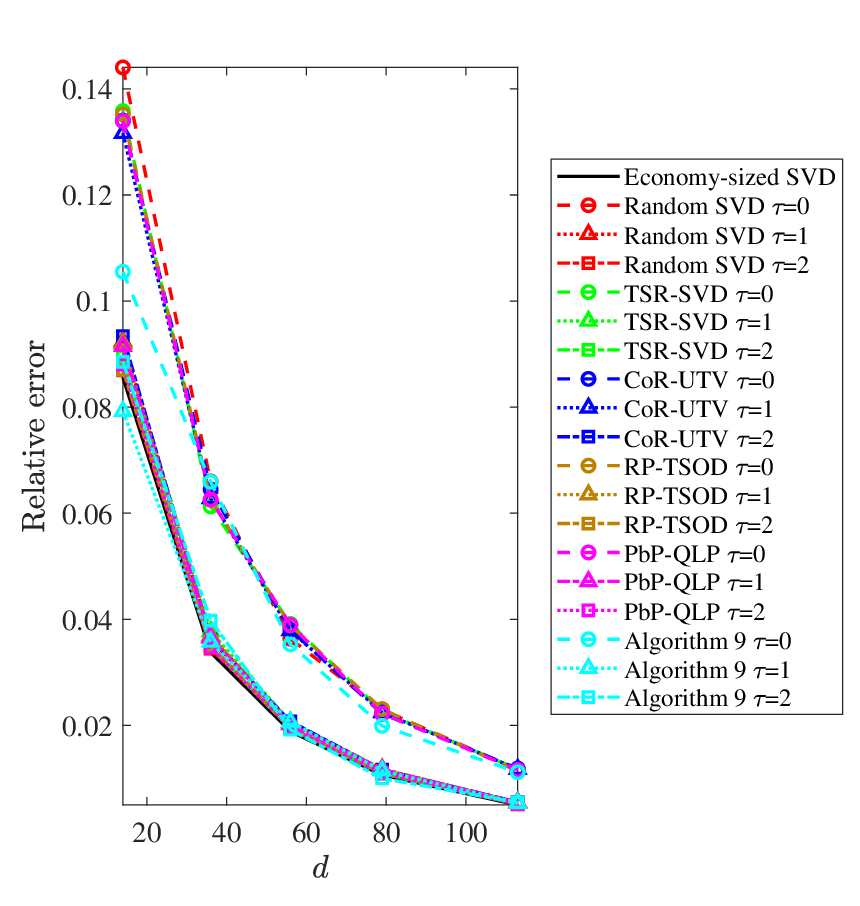}\\
  \caption{Low-rank image reconstruction Frobenius norm approximation relative error. This
figure displays the relative errors incurred by different methods in reconstructing the
house image.}\label{example-house-err}
\end{figure}

\begin{figure}[!ht]
  \centering
  \includegraphics[width=3.3in]{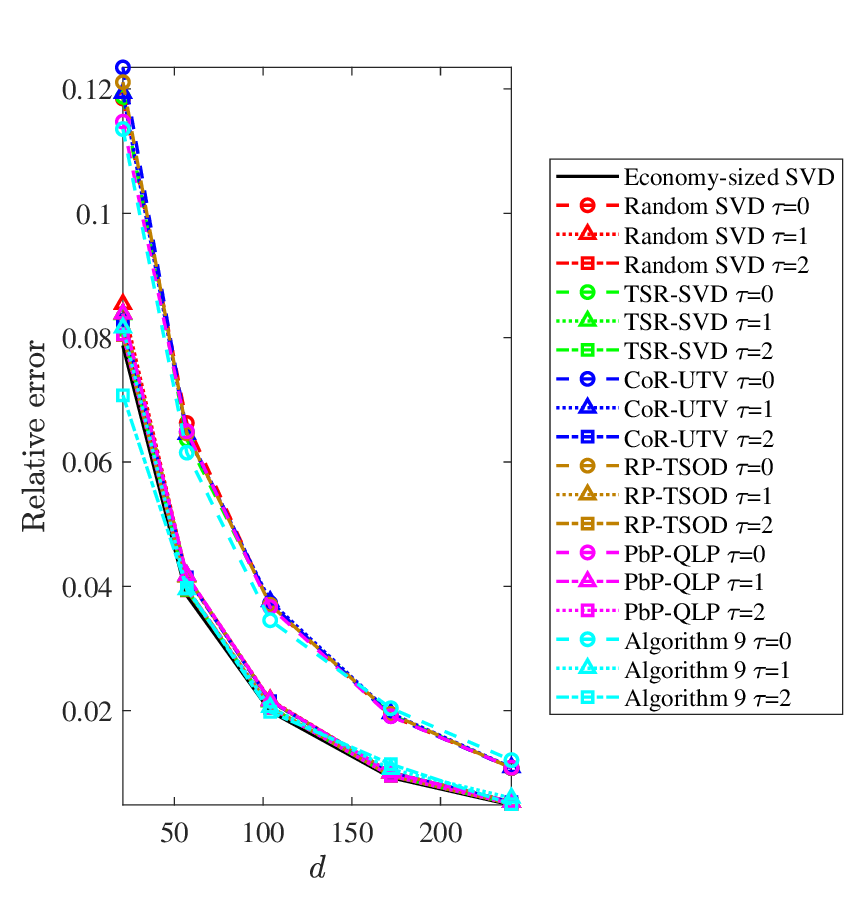}\\
  \caption{Low-rank image reconstruction Frobenius norm approximation relative error. This
figure displays the relative errors incurred by different methods in reconstructing the
jetplane image.}\label{example-jetplane-err}
\end{figure}

The goal of this experiment is to assess the performances of Algorithm \ref{alg:A} on real world data.
The real image data set is adopted, which is available at \url{https://github.com/Whisper329/image-processing}.
This dataset, extensively employed in image processing, comprises 12 grayscale images each of 512x512 resolution.
Our aim is to reconstruct these low-rank gray images.
Furthermore, we compare our results against those constructed by several other algorithms from the literature.
The specific operation is as follows: Since the calculation results of Algorithm \ref{alg:A} are closely related to the given precision $\varepsilon$, during the experiment, we first fixed the $\varepsilon$ value to 0.001 for all images, thereby obtaining the value of the sampling parameter $d$ of Algorithm \ref{alg:A}, and then substitute random SVD, TSR-SVD, CoR-UTV, RP-TSOD, PbP-QLP and the built-in function economy-sized SVD in MATLAB. In this way, each image selects different sampling parameters for reconstruction, and the reconstruction results are shown in Figures \ref{example-house} and \ref{example-jetplane}.
Figures \ref{example-house} and \ref{example-jetplane} show the original image and the reconstructed image.
We reconstructed these 12 gray images 1000 times using 7 methods respectively. The calculation time is shown in Table \ref{example-t}.
Figures \ref{example-house-err} and \ref{example-jetplane-err} show the relative error changes of images reconstructed by image house.tif for different sample sampling parameters $d$.

We can draw the following conclusions from the Figures \ref{example-house}, \ref{example-jetplane}, \ref{example-house-err}, \ref{example-jetplane-err} and Table \ref{example-t}.
\begin{itemize}
  \item The approximation of Algorithm \ref{alg:A} is visually as good as that of the algorithms mentioned earlier. All methods can well restore the feature information of the original image.
      It can be seen from Figures \ref{example-house} and \ref{example-jetplane} that as the power factor $\tau$ is selected larger, the reconstructed image becomes clearer and the image quality is closer to the original image.
  \item For the same power factor $\tau$, the computational time of reconstructing images through Algorithm \ref{alg:A} is the fastest compared with other algorithms, several times that of the second fastest algorithm.
  \item Figures \ref{example-house-err} and \ref{example-jetplane-err} show the relative errors of different methods when reconstructing the house.tif image. It can be observed that, for the same power factor $\tau$, the relative errors of all mentioned algorithms are almost the same, and as the sampling parameter increases, the relative error values become smaller and smaller.
      In addition, by increasing the power factor $\tau$, the relative errors of Algorithm \ref{alg:A} is close to that of the built-in function economy-sized SVD in MATLAB.
\end{itemize}


\begin{table}[!ht]
  \centering
  \caption{The calculation time of different methods}\label{example-t}
  \begin{tabular}{lcc}
    \hline
    Algorithms  & & Calculation time (second) \\
    \hline
    Economy-sized SVD                  &          & 38.34\\
                                       &          &       \\
    \multirow{3}{*}{Random SVD}        & $\tau=0$ & 8.85\\
                                       & $\tau=1$ & 12.03\\
                                       & $\tau=2$ & 15.39\\
                                       &          &       \\ 
    \multirow{3}{*}{TSR-SVD}           & $\tau=0$ & 9.84\\
                                       & $\tau=1$ & 13.18\\
                                       & $\tau=2$ & 16.41\\
                                       &          &       \\ 
    \multirow{3}{*}{CoR-UTV}           & $\tau=0$ & 10.41\\
                                       & $\tau=1$ & 13.66\\
                                       & $\tau=2$ & 16.89\\
                                       &          &       \\ 
    \multirow{3}{*}{RP-TSOD}           & $\tau=0$ & 9.05\\
                                       & $\tau=1$ & 12.41\\
                                       & $\tau=2$ & 15.78\\
                                       &          &       \\ 
    \multirow{3}{*}{PbP-QLP}           & $\tau=0$ & 11.43\\
                                       & $\tau=1$ & 14.90\\
                                       & $\tau=2$ & 18.07\\
                                       &          &       \\ 
    \multirow{3}{*}{Algorithm \ref{alg:A}} & $\tau=0$ & \textbf{5.94} \\
                                       & $\tau=1$ & 9.09\\
                                       & $\tau=2$ & 12.64\\
    \hline
  \end{tabular}
\end{table}

\section{Conclusion}
We propose fast low-rank approximation algorithm for unknown rank. The introduced randomized algorithms for basis extraction allows the acquisition of the matrix rank, addressing the limitations of empirically obtaining numerical ranks. The presented algorithms undergo theoretical analysis, and numerical experiments on random matrices indicate that the proposed algorithm exhibits superior time efficiency while maintaining accuracy comparable to existing algorithms.  Power method techniques are incorporated into the algorithms to enhance precision, and experimental results validate the effectiveness of this enhancement.
Finally, the proposed algorithm is well applied to image reconstruction, yielding favorable results and outperforming the second fastest algorithm by approximately several times in terms of speed.


\end{document}